\documentclass[10pt,a4paper]{article}
\usepackage{amsfonts,dsfont,amsmath,amssymb,amsthm,enumerate,color}
\usepackage{graphicx}
\usepackage{fullpage}
\usepackage{authblk}
\usepackage{bold-extra}

\definecolor{orange}{rgb}{1,0.25,0}
\definecolor{darkgreen}{rgb}{0.0, .5, 0.13}

\usepackage[colorlinks=true,linkcolor=blue,citecolor=red, pdfborder={0 0 0}]{hyperref}

\renewcommand{\t}{\theta}
\renewcommand{\i}{\mathrm{i}}

\newcommand{\ep}{\varepsilon}

\newcommand{\CC}{\mathbb{C}}
\newcommand{\NN}{\mathbb{N}}
\newcommand{\RR}{\mathbb{R}}

\newcommand{\E}{\mathcal{E}}
\renewcommand{\L}{\mathcal{L}}
\renewcommand{\H}{\mathcal{H}}

\newcommand{\N}{\mathcal{N}}

\renewcommand{\P}{\mathcal{P}}

\newcommand{\vertiii}[1]{{\left\vert\kern-0.25ex\left\vert\kern-0.25ex\left\vert #1 
		\right\vert\kern-0.25ex\right\vert\kern-0.25ex\right\vert}}

\renewcommand{\i}{{\rm i}}

\theoremstyle{plain}
\newtheorem{theorem}{Theorem}[section]
\newtheorem{lemma}[theorem]{Lemma}
\newtheorem{proposition}[theorem]{Proposition}
\newtheorem{corollary}[theorem]{Corollary}

\theoremstyle{definition}

\theoremstyle{plain}
\newtheorem{remark}[theorem]{Remark}

\numberwithin{equation}{section}

\title{\vspace{-1.8cm} Title
}

\renewcommand{\t}{\theta}

\begin{document}

\title{Quantitative high-contrast homogenisation of composite media with highly anisotropic fibres}
\author{Shane Cooper and Ilia Kamotski} 

\maketitle

\begin{abstract}
This work concerns a non-uniformly elliptic second-order differential operator with periodic coefficients that describes a composite medium formed by a periodic array of highly anisotropic cylindrical fibres set in an isotropic matrix. The ratio of transverse to axial conductivity in the fibres is linked to the period of the coefficients through a `critical' high-contrast scaling, and the axial conductivity is allowed to vary along the fibres. As the period $\varepsilon$ tends to zero, ellipticity degenerates in the directions transverse to the fibres. We are concerned with the resolvent of this operator in the small-$\varepsilon$ limit. For fibres with constant axial conductivity, the two-scale limit of such problems was identified by Cherednichenko, Smyshlyaev and Zhikov (Proceedings of the Royal Society of Edinburgh: Section A Mathematics 136(1), 87--114 (2006)), and a full two-scale asymptotic expansion, with remainder estimates for smooth right-hand sides, was constructed by Cherednichenko (Asymptotic Analysis 49(1--2), 39--59 (2006)). In the non-constant axial conductivity setting, considered here, we approximate the resolvent in terms of the two-scale limit operator and prove that the error is order $\varepsilon$ in the operator norm, uniformly over right-hand sides in the unit ball of $L^2(\mathbb{R}^3)$; we also obtain order-$\varepsilon$ estimates for the reciprocals of the spectral band functions. The proof follows the general strategy of Cooper, Kamotski and Smyshlyaev (preprint available at https://arxiv.org/abs/2307.13151). Two features of the present problem, however, fall outside that strategy: the loss of ellipticity depends on direction, and a key `spectral gap' assumption of that work is satisfied only in a weaker form.
Additionally, for non-constant axial conductivity,
a boundary-layer analysis at the surface of each fibre is required to recover order-$\varepsilon$ estimates. This work builds on results announced, without proof, for fibres with
circular cross-section in Cooper and Kamotski (Direct and Inverse Problems with Applications: Extended Abstracts of the 2024 GAP Center Summer School, pp.~69--75, Birkh\"auser (2026)).
\end{abstract}

\section{Introduction}

In this work, we analyse the asymptotic behaviour, as \(\varepsilon\to0\), of
the solution \(u_\varepsilon\in H^1(\mathbb R^3)\) to the problem
\begin{equation}\label{p.strong}
\int_{\RR^3}A_\varepsilon\big(\tfrac{\cdot}{\varepsilon}\big)\nabla u_\varepsilon \cdot \overline{\nabla \phi}
     +\int_{\RR^3}u_\varepsilon \overline{\phi} =\int_{\RR^3}f\overline{\phi} ,\qquad \forall \phi \in H^1(\RR^3),
\end{equation}
for given  $f\in L^2(\mathbb R^3)$.

Here, as is typical in homogenisation problems, we take the coefficient matrix $A_\ep$ to be $\square$-periodic, with $\square:=[-\frac{1}{2},\frac{1}{2})^3$, so that $A_\ep(\frac{x}{\ep})$ models the material properties of a periodic anisotropic composite medium with small period $\ep$. In addition, the coefficients degenerate with the critical scaling familiar from now-classical `high-contrast' homogenisation problems (cf. \cite{ADH}, \cite{All}, \cite{Zh}); here, however, we have a non-classical anisotropic degeneracy. Namely, the composite consists of highly anisotropic cylindrical fibres embedded in a background medium: the transverse `conductivity' of the fibres is of order $\ep^2$, whereas the axial conductivity of the fibres, and the conductivity of the background medium, is of order one.  Such operators serve as a prototype for elastic composites whose fibres are stiff along their axis and soft across it and as a model for wave propagation and localisation in highly anisotropic periodic composites, cf.~\cite{VPS}.
Precisely, the periodic function $A_\ep$ is defined on its periodic cell $\square$ as follows:
\[
    A_\varepsilon(y)=\big(1-\chi_0(y')\big)I+
     \chi_0(y')\operatorname{diag}\bigl(\varepsilon^2,\varepsilon^2,a(y_3)\bigr), \quad y=(y',y_3),\, y' \in \square' = [-\tfrac12,\tfrac12)^2,\, y_3 \in [-\tfrac{1}{2},\tfrac{1}{2}),
\]
for \(\chi_0\) the characteristic function of \(B\), where \(B\) is a
simply connected domain with smooth boundary whose closure is contained in $(-\tfrac12,\tfrac12)^2$.  The coefficient \(a\) is measurable, \(1\)-periodic,
and satisfies \(\nu\leq a\leq\nu^{-1}\) almost everywhere for some \(\nu>0\).


Problems of this type have previously been studied.
\cite{ChSmZh} studied this problem on a
macroscopic periodic cell, for a constant axial conductivity $a$.\footnote{In  \cite{ChSmZh} it was commented that the constant axial conductivity assumption was  not essential to their approach and could be removed.} They established weak
two-scale convergence (see for example \cite{All}, \cite{Ng}, \cite{Zh}) to a limit of the form \(u(x)+v(x,y')\) and showed that
eliminating the fibre component \(v\) gives an equation for \(u\) which is
non-local in the direction of the fibres. Subsequently, \cite{ChAsy}  constructed a full two-scale
asymptotic expansion, for constant $a$ and smooth periodic right-hand sides,  and obtained estimates for its remainder with constants that (implicitly) depend on the right-hand-side.

Our aim is to establish, for non-constant $a$, order-$\ep$ resolvent estimates that are uniform with respect to the right-hand side in the unit ball of $L^2(\RR^3)$.
 Such operator-norm estimates are by now well developed in classical (uniformly elliptic) homogenisation, see e.g.\ \cite{BS04,ZhP05,KLS12}, and have more recently been obtained for high-contrast problems, see e.g.\ \cite{ChCo16,ChErKi20,CoKaSm}. We will also obtain order-$\ep$  asymptotics for the reciprocals of the spectral band functions of the operator $L_\ep$  generated by problem \eqref{p.strong}, uniform in the quasi-momenta and period $\ep$.

To achieve this goal, we extend the work \cite{CoKaSm} to this setting. Such an extension presents two main technical difficulties which we shall now describe. First we reformulate \eqref{p.strong} in the usual way. By applying the
unitary rescaling \((\Gamma_\varepsilon f)(x)=\varepsilon^{3/2}f(\varepsilon x)\)
and the Floquet-Bloch-Gelfand transform $U$,\footnote{ $U$ is given by the continuous extension of the mapping 
\[
  (Uf)(\theta,y)=(2\pi)^{-3/2}
       \sum_{m\in\mathbb Z^3} f(y+m)e^{-i\theta\cdot(y+m)},
  \qquad \theta\in [-\pi,\pi)^3.
\]
}
 we
obtain, for almost every $\t \in \square^* := [-\pi,\pi)^3$, the problem
\begin{equation}\label{p0} \left\{
 \begin{aligned}&\text{Find $
 u_{\varepsilon,\theta}\in H^1_{\rm per}(\square)$ the solution to}\\
 &\varepsilon^{-2}a_\theta(u_{\varepsilon,\theta},\widetilde u)
   +b_\theta(u_{\varepsilon,\theta},\widetilde u)
       =(f_{\varepsilon,\theta},\widetilde u),
 \qquad \forall \widetilde u\in H^1_{\rm per}(\square).
 \end{aligned} \right.
\end{equation}
Here,  $f_{\varepsilon,\theta}=(U\Gamma_\varepsilon f)(\theta,\cdot) \in L^2(\square)$, $(\cdot,\cdot)$ denotes the inner product of $L^2(\square)$ and the forms are given by
\begin{align}
 a_\t(u,u)=:a_\theta[u]&=\int_{\square\setminus F}|(\nabla+i\theta)u|^2\,dy
             +\int_F a(y_3)|(\partial_3+i\theta_3)u|^2\,dy,\label{at}\\
b_\t(u,u)=: b_\theta[u]&=\int_F|(\nabla'+i\theta')u|^2\,dy
             +\int_\square |u|^2\,dy,\nonumber
\end{align}
where $F= B \times(-\tfrac{1}{2},\tfrac{1}{2})$, $\nabla'=(\partial_1,\partial_2)$ and $\theta'=(\theta_1,\theta_2)$.
Henceforth, we will consider the problem \eqref{p0} for arbitrary right-hand-side $f$ belonging to $L^2(\square)$.

The first difficulty in  applying the approach of \cite{CoKaSm} to
\eqref{p0} is related to a `spectral-gap condition'  which states that, for each fixed
\(\theta\), the form \(a_\theta\)  is coercive in $H^1_{per}(\square)$ (after projecting away from its kernel).  Amongst other things, this particular condition ensured that first order correction of $u_{\ep,\t}$ lies in $H^1_{per}(\square)$ which is crucial to the analysis in \cite{CoKaSm}.  This assumption, which holds for many examples of interest, cf. \cite{IKVPS} and \cite{CoKaSm}, fails in the present example, cf. \cite{IKVPS}. However, a weaker spectral-gap condition holds.
Indeed,  consider
\begin{equation}\label{V}
V  
 = 	\{ v \in H^1_{per}(\square) : \text{$v$ constant in $\square \backslash F$ and  $v$ independent of $y_3$ in $F$} \} = \CC \dot{+} H^1_0(B),
\end{equation}
where $H^1_0(B)$ is understood as a subset of $H^1_{per}(\square)$ in the obvious way.
Now, denote  the closed subset $W $ of $H^1_{per}(\square)$ that is the orthogonal complement of $V$  in some $H^1$ equivalent inner product (whose exact definition will be given later). Then, it can be shown that 
 $a_\t$ satisfies 
the following weak coercivity estimate:  there exists a positive constant $\gamma$ such that
	\begin{align*}
	a_\t[w] &\ge  \gamma \|w\|^2_{L^2(\square)}, \quad \forall w \in W, \ \forall \t \in \square^*. 
\end{align*}
Overcoming this absence of full coercivity is one of the main challenges of this current work as it leads to  new analysis not present in \cite{CoKaSm}; in particular, we need to investigate internal boundary-layer effects near each fibre.
The second difficulty is the directional dependence of the degeneracy. For $a_\t$, given by \eqref{at}, one has
\[
 a_\theta[u]\gtrsim
 |\theta|^2\|u\|_{L^2(\square\setminus F)}^2
       +|\theta_3|^2\|u\|_{L^2(\square)}^2 .
\]
This directional dependence on the degeneracy of the form $a_\t$ is not present in typical homogenisation problems, though it arises in other high-contrast settings, cf. \cite{Co}. Such degeneracies are not covered in the general framework of  \cite{CoKaSm}.

To state our results, let \(\mathcal H\) be the closure of \(V\) in
\(L^2(\square)\) and let \(\mathcal P:L^2(\square)\to\mathcal H\) be the
orthogonal projection. For \(\xi\in\mathbb R^3\), denote by \(\mathcal L_\xi\)
the operator in \(\mathcal H\) associated with the form
\[
\begin{aligned}
A_\xi[c+z]={}&(A^{\rm h}\xi\cdot\xi)|c|^2
       +a^{\rm h}|\xi_3|^2\|c+z\|_{L^2(B)}^2 +\|\nabla' z\|_{L^2(B)}^2+\|c+z\|_{L^2(\square)}^2,
       \ c+z\in V=\CC \dot{+} H^1_0(B) .
\end{aligned}
\]
The operators \(\mathcal L_\xi\) are the Fourier fibres, in the
macroscopic variable, of the two-scale operator \(\mathbb L\) defined (precisely in Section \ref{sec.twoscale}) by the problem: for given $g(x,y') \in L^2(\RR^3;\H) \subset  L^2(\RR^3 \times \square)$, find $u(x) + v(x,y')$ the solution to
\begin{equation}\label{p4}
\begin{cases}
-\operatorname{div}_x(A^{\rm h}\nabla_xu)
 -a^{\rm h}\partial_{x_3}^{\,2}\bigl(|B|u+\langle v\rangle\bigr)
 +u+\langle v\rangle=\langle g \rangle, & x\in\mathbb R^3,\\
-\Delta_{y'}v-a^{\rm h}\partial_{x_3}^{\,2}(u+v)+u+v=g,
   & x\in\mathbb R^3,\ y'\in B,\\
v(x,y')=0, & x\in\mathbb R^3,\ y'\in\partial B.
\end{cases}
\end{equation}
Here, $\langle h\rangle(x):=\int_{\square'} h(x,y')\,dy'$, $A^h$ and $a^h$ are the homogenised coefficients defined in Section \ref{sec.limform}.  Problem \eqref{p4} is the two-scale problem associated with \eqref{p.strong}  (cf. \cite{ChSmZh}). Additionally, let \(L_{\varepsilon,\theta}\) be the operator associated with \eqref{p0},  
and let \(\mathcal E_{\theta'}\) denote multiplication by
\(e^{i\theta'\cdot y'}\). Our first main result,  Theorem~\ref{thm1.2}, states that 
\[
 \bigl\|L_{\varepsilon,\theta}^{-1}
 -\mathcal E_{\theta'}^*\mathcal P^*\mathcal L_{\theta/\varepsilon}^{-1}
       \mathcal P\mathcal E_{\theta'}\bigr\|_{L^2(\square)\to L^2(\square)}
 \le C\varepsilon,\qquad \theta\in\square^*,\quad 0<\varepsilon<1 .
\]
Theorem \ref{thm1.3}, compares spectral band functions. If \(\lambda_\varepsilon^{(k)}(\theta)\) and
\(\Lambda^{(k)}(\xi)\) are the eigenvalues of \(L_{\varepsilon,\theta}\) and
\(\mathcal L_\xi\), in increasing order and repeated according to
multiplicity, then one has
\[
 \left|\frac{1}{\lambda_\varepsilon^{(k)}(\theta)}
      -\frac{1}{\Lambda^{(k)}(\theta/\varepsilon)}\right|
 \le C\varepsilon,\qquad k\in\mathbb N,\quad\theta\in\square^* .
\]
Our final result, in Section \ref{sec.twoscale}, relates $L_\ep^{-1}$ (defined on $L^2(\RR^3)$) to $\mathbb{L}^{-1}$ (defined on a closed subspace  of $L^2(\RR^3\times \square)$) with the help of the isometry (see \cite{CoKaSm} and, for a similar construction without $\E$, see \cite{Well})
\[
 J_\varepsilon:=\Gamma_\varepsilon^{-1}\mathcal F^{-1}
       \chi\mathcal E U\Gamma_\varepsilon
 \colon L^2(\mathbb R^3)\longrightarrow L^2(\mathbb R^3\times\square),
\]
where \(\mathcal E\) acts fibrewise as \(\mathcal E_{\theta'}\), \(\chi\)
extends by zero outside \(\square^*\), and \(\mathcal F\) is the Fourier
transform in the first variable. Upon inverting the transforms, $U$ and $\Gamma_\ep$, we obtain, in Theorem~\ref{thm2}, the
whole-space estimate 
\[
 \bigl\|L_\varepsilon^{-1}
      - J_\varepsilon^*P^*\mathbb L^{-1}P J_\varepsilon
 \bigr\|_{L^2(\mathbb R^3)\to L^2(\mathbb R^3)}
 \le C\varepsilon,\qquad 0<\varepsilon<1 .
\]
 Here, \(P\) is the projection acting as \(\mathcal P\) in the
cell variable.

Some of the results of this work, for the case of $B$ being a disc, were announced (without proof) in \cite{CoKaPRE}.

Throughout the article, the notation
\[
X \lesssim Y,
\]
means
\[
X \le C Y
\]
for some constant $C>0$ independent of $X$ and $Y$ and the parameters $\ep$, $\t\in \square^*$ .
\section{A priori estimates}
We begin with some inequalities for the forms $a_\t$ and $b_\t$, uniform in $\t\in\square^*$, which will be used throughout the article.
\begin{proposition}\label{prop.directionaldep}
\begin{equation*} \begin{aligned}
     |\t_3|^2 \int_{\square}|u|^2  + |\t|^2 &\int_{\square\backslash F}|u|^2  + \int_{\square\backslash F}|\nabla u|^2 + \int_{\square}  |\partial_3u|^2 \\
     &\lesssim \int_{\square \backslash F} |(\nabla+ \i \t)u|^2 +\int_{\square}  |(\partial_3+ \i \t_3)u|^2, \quad \forall u \in H^1_{per}(\square), \forall \t \in \square^*. 
     \end{aligned}
\end{equation*}
\end{proposition}
\begin{proof}
    A standard Fourier series argument gives
    \begin{equation*}
  |\t_3|^2 \int_{\square} | u |^2 \le \int_{\square} |(\partial_3 + \i \t_3 ) u |^2\quad \forall u \in H^1_{per}(\square). 
    \end{equation*}
    Additionally, one has (see for example  \cite[Equation (7.42)]{CoKaSm}) that
    \begin{equation}\label{27.3.26e3}
|\t'|^2 \int_{\square'\backslash B} | \phi |^2 \lesssim \int_{\square'\backslash B} |(\nabla' + \i \t' ) \phi |^2 , \quad \forall \phi \in H^1_{per}(\square'\backslash B).
    \end{equation}
    These two facts, after integrating \eqref{27.3.26e3} with respect to $y_3$, imply the desired result.
\end{proof}
\begin{proposition}\label{Prop127.3.26}
For any $u\in H^1_{per}(\square')$ that satisfies $\int_{\square' \backslash B} u =0$  and 
    \begin{equation}\label{27.3.26e1}
\int_{B}  \nabla'u \cdot \overline{\nabla' v} = 0, \quad \forall v \in H^1_0(B),
    \end{equation}
    one has 
    \begin{equation}\label{27.3.26e4}
  \int_{\square'}| u|^2 +\int_{\square'}| (\nabla' + \i \t' )u|^2 \lesssim \int_{\square' \backslash B}| (\nabla' + \i \t' )u|^2, \quad \forall \t' \in [-\pi,\pi)^2.
    \end{equation}
\end{proposition}
\begin{proof}
  First note that, one has
   \begin{equation}\label{27.3.26e2}
 \int_{\square'}| u|^2 +\int_{\square'}| \nabla' u|^2 \lesssim \int_{\square' \backslash B}| \nabla' u|^2 + \int_{\square'\backslash B}| u|^2.
   \end{equation}

   This readily follows, for example, by taking test function $v= u - Eu$ in \eqref{27.3.26e1} where $E:H^1(\square'\backslash B) \rightarrow H^1(\square')$ is a Sobolev extension. 

By \eqref{27.3.26e2} and the Poincar\'{e} inequality for $\square' \backslash B$, it follows that \eqref{27.3.26e4} holds for $\t'=0$; namely,
one has
 \begin{equation*}
 \int_{\square'}| u|^2 +\int_{\square'}| \nabla' u|^2 \lesssim \int_{\square' \backslash B}| \nabla' u|^2 .
   \end{equation*}
Therefore,
 \begin{equation*}
 \int_{\square'}| u|^2 +\int_{\square'}| (\nabla'+\i \t') u|^2 \lesssim \int_{\square' \backslash B}| \nabla' u|^2 , \quad \forall\, \t' \in [-\pi, \pi)^2.  
   \end{equation*}  
It remains to notice that
\begin{equation*}
\int_{\square'\backslash B}| \nabla' u|^2 \le 2 \left(\int_{\square' \backslash B}| (\nabla' + \i \t' )u|^2 + |\t'|^2\int_{\square' \backslash B}|  u|^2\right),
\end{equation*}
and then, upon recalling \eqref{27.3.26e3},  the desired result, inequality \eqref{27.3.26e4}, follows.
\end{proof}

Recall the closed linear space $V$, given by \eqref{V}, and $W$ the orthogonal complement of $V$ in $H^1_{per}(\square)$ with respect to the equivalent inner product
\begin{equation}\label{def.W}
( \tilde{u},\tilde{v} )_{H^1(\square)} : = \int_{\square \backslash F} \tilde{u} \overline{\tilde{v}} + \int_{\square}  \nabla \tilde{u}  \cdot \overline{\nabla \tilde{v}}.
\end{equation}
Notice that
\begin{equation}\label{W}
W = \left\{ w \in H^1_{per}(\square) \, : \, \text{$\int_{\square \backslash F} w =0$ and $\Delta' \langle w\rangle_3 = 0$ in $B$}    \right\},
\end{equation}
where we define $\langle f \rangle_3(y') : = \int_{-\tfrac{1}{2}}^{\tfrac{1}{2}} f(y',y_3) \, \mathrm{d}y_3$.  

Our first observation is that $a_\t$ is uniformly coercive on $W$ in the following sense:
\begin{proposition}\label{prop.Wuniformcoercive}
    \[
 \|w\|_{L^2(\square)}^2  + \|\partial_3 w\|_{L^2(\square)}^2+ \| \langle w \rangle_3\|^2_{H^1(\square')} + \|\nabla w\|_{L^2(\square \backslash F)}^2 \lesssim a_\t[w] , \quad \forall w \in W,\, \forall \t \in \square^*.
    \]
\end{proposition}
 \begin{proof}
For ease of exposition, let us denote 
\[
Q_\t[w]  : =   		\int_{\square\backslash F} | (\nabla + \i \t) w(y)|^2\, \mathrm{d}y + \int_\square | (\partial_3 +\i \theta_3) w(y)|^2 \, \mathrm{d}y.
\]
Clearly, since $a$ is positive, $Q_\t$ is uniformly bounded by $a_\t$. Now, Proposition \ref{prop.directionaldep} implies
\[
 \|\partial_3 w\|_{L^2(\square)}^2 + \|\nabla w\|_{L^2(\square \backslash F)}^2 \lesssim Q_\t[w].
\]
Therefore, it remains to show
\begin{equation}\label{15.4.26e1}
\|w\|_{L^2(\square)}^2 + \|\langle w \rangle_3 \|_{H^1(\square')}^2
\lesssim Q_\t[w] . 
\end{equation}
The decomposition $w(y) = w_1(y') + w_2(y)$ where $\langle w_2 \rangle_3=0$ holds. Notice that
\[
Q_\t[w] = Q_\t[w_1]+Q_\t[w_2], \quad \text{and} \quad \|w\|_{L^2(\square)}^2=\| w_1 \|_{L^2(\square')}^2+\|w_2\|_{L^2(\square)}^2.
\]
Clearly, one has
\[
\pi^2 \|w_2\|_{L^2(\square)}^2 \le Q_\t[w_2] .
\]
Furthermore, since $w \in W$ then the assumptions of  Proposition \ref{Prop127.3.26} are satisfied for $u=w_1 = \langle w \rangle_3$. Consequently, it follows from \eqref{27.3.26e4}, and the fact $\t \in [-\pi,\pi)^3$, that
\[
\|  w_1 \|_{H^1(\square')}^2  \lesssim Q_\t[w_1].
\]
Hence, \eqref{15.4.26e1} holds.
\end{proof}

\begin{proposition}\label{prop.bcrossterm}
\[
    \left| b_\t(v,w)  \right| \lesssim \sqrt{b_\t[v]} \sqrt{a_\t[w]}, \quad \forall v \in V, \forall w \in W, \forall \t \in \square^*.
\]
\end{proposition}
\begin{proof} 
Since $v$ is independent of $y_3$, and by Proposition \ref{prop.Wuniformcoercive}, one has
\begin{flalign*}
    | b_\t(v,w) | & = |  b_\t(v,\langle w\rangle_3)|  \le \sqrt{b_\t[v]}\sqrt{b_\t[\langle w \rangle_3]}\lesssim\sqrt{b_\t[v]}\sqrt{a_\t[ w ]}.
\end{flalign*}
\end{proof}

\begin{proposition}\label{prop.weakspectralgap}
   \[
a_\t[v]+a_\t[w]\lesssim a_\t[v+ w], \quad \forall v \in V, \forall w \in W, \forall \t \in \square^*.
\]
\end{proposition}
\begin{proof}Clearly it is enough to prove that
\[
a_\t[v] \lesssim a_\t[v+w].
\]
Notice that
  \begin{equation}\label{eqnatv1}
    a_\t[v] =  |\t|^2\int_{\square \backslash F} |v|^2 +   |\t_3|^2\int_{ F} a(y_3) |v|^2.        
    \end{equation}
By Proposition \ref{prop.directionaldep}, one has
    \begin{flalign}
 a_\t[v + w] &\gtrsim \| \nabla (v+w)\|_{L^2(\square\backslash F)}^2+ |\t|^2 \|v + w\|_{L^2(\square\backslash F)}^2 \nonumber \\
 &=  \| \nabla w\|_{L^2(\square\backslash F)}^2+  |\t|^2 \| v\|_{L^2(\square\backslash F)}^2  + |\t|^2 \|  w\|_{L^2(\square\backslash F)}^2, \label{3.10.26e1}
 \end{flalign}
where the equality comes from the fact that 
 $v$ is constant in $\square \backslash F$ and $\int_{\square \backslash F} w =0$.

 On the other hand, Proposition \ref{prop.directionaldep} implies
 \begin{flalign}
a_\t[v + w] & \gtrsim  |\t_3|^2 \|v + w\|_{L^2(F)}^2 \nonumber \\   & \ge \frac{1}{2} |\t_3|^2 \| v \|_{L^2(F)}^2 - 2|\t_3|^2 \| \langle w\rangle_3\|_{L^2(F)}^2, \label{3.10.26e2}      
 \end{flalign}
where in the last inequality we used the fact that $v$ is independent of $y_3$.

Finally,  \eqref{27.3.26e2} (for $u = \langle w \rangle_3$) implies that
\[
\int_{\square'}| \langle w \rangle_3|^2 +\int_{\square'}| \nabla' \langle w \rangle_3|^2 \lesssim \int_{\square' \backslash B}| \nabla' \langle w \rangle_3|^2 + \int_{\square'\backslash B}| \langle w \rangle_3|^2,
\]
which readily gives
\[
 \| \langle w\rangle_3\|_{L^2(F)}^2  \lesssim \| \nabla w\|_{L^2(\square\backslash F)}^2 + \| w\|_{L^2(\square\backslash F)}^2.
\]
This inequality, along with \eqref{3.10.26e2} and $\t_3 \in [-\pi,\pi)$, gives
\begin{equation}\label{3.10.26e3}
a_\t[v + w]  + \| \nabla w \|_{L^2(\square\backslash F)}^2 +|\t_3|^2 \| w \|_{L^2(\square\backslash F)}^2\gtrsim |\t_3|^2 \| v \|_{L^2(F)}^2.
\end{equation}
Summing up \eqref{3.10.26e1} and \eqref{3.10.26e3}, after multiplying by the appropriate constants, gives
\[
a_\t[v + w]  \gtrsim|\t|^2\int_{\square \backslash F} |v|^2 +   |\t_3|^2\int_{ F} |v|^2,
\]
which along with \eqref{eqnatv1} proves the result.
\end{proof}
\begin{corollary} \label{corAonV}
     \[
   |\t|^2 \| v\|_{L^2(\square\backslash F)}^2 + |\t_3|^2 \| v\|_{L^2(F)}^2 \lesssim a_\t[v+ w], \quad \forall v \in V, \forall w \in W, \forall \t \in \square^*.
\]
\end{corollary}

\begin{proposition}\label{prop.bsplit}
\[
b_\t[v]+b_\t[w] \lesssim b_\t[v+w], \quad \forall v \in V, \forall w \in W, \forall \t \in \square^*.
\]
\end{proposition}
\begin{proof}Clearly, it is enough to prove
\[
b_\t[v] \lesssim b_\t[v+w], \quad \forall v \in V, \forall w \in W, \forall \t \in \square^*.
\]
    Since 
    \begin{equation}\label{btequiv}
b_0[u] \lesssim b_\t[u] \lesssim b_0[u], \quad \forall u \in H^1(\square), \forall \t \in \square^*, 
    \end{equation}
 it suffices to prove the desired result for $\t=0$.

To that end, notice that
\begin{flalign}
b_0[v+w] & \ge \int_F | \nabla' v + \nabla' w|^2 + \int_{\square \backslash F} |v+w|^2 \nonumber \\
& = \int_F | \nabla' v|^2 +  \int_F | \nabla' w|^2+\int_{\square \backslash F} |v|^2  + \int_{\square \backslash F} |w|^2 \label{16.4.26e1}
\end{flalign}
where the last equality used the definitions of $V$ and $W$ (cf. \eqref{V} and \eqref{def.W}). 
Now, since $v = v_0 +c$ for $v_0 \in H^1_0(B)$ and $c \in \CC$, then
\begin{flalign*}
    \| v\|_{L^2(F)}^2 &  \le 2| F| |c|^2 +2 \| v_0 \|_{L^2(F)}^2 \\
&    \lesssim |\square \backslash F| |c|^2 + \| v_0 \|_{L^2(F)}^2 \\
    & =  \| v \|_{L^2(\square \backslash F)}^2+  \| v_0 \|_{L^2(F)}^2 \\
   & \lesssim   \| v \|_{L^2(\square \backslash F)}^2+  \| \nabla' v_0 \|_{L^2(F)}^2 \\
   & =  \| v \|_{L^2(\square \backslash F)}^2+  \| \nabla' v\|_{L^2(F)}^2
\end{flalign*}
which along with \eqref{16.4.26e1} gives
\[
b_0[v] = \int_F | \nabla' v|^2 + \int_{\square} |v|^2 \lesssim  b_0[v+w].
\]
\end{proof}
For the form 
\begin{equation}\label{A}
    A_{\ep,\t} : = \ep^{-2} a_\t + b_\t,
\end{equation}
one has 
\begin{equation}\label{cor1}
A_{\ep,\t}[v]  + A_{\ep,\t}[w]  \lesssim A_{\ep,\t}[v+ w], \quad \forall v \in V, \forall w \in W,\, \forall \t \in \square^*.
\end{equation}

\section{Elimination of the coercive subspace $W$}
We introduce the bounded linear operator $\mathcal{N}_{\ep,\theta}: V \rightarrow W$, defined by the mapping $v \mapsto \mathcal{N}_{\ep,\theta}v$ the solution to:
\begin{equation}\label{firstcorrector}
\ep^{-2} a_\t(\mathcal{N}_{\ep,\theta} v, \tilde{w}) + b_\t(\mathcal{N}_{\ep,\theta} v, \tilde{w})    =  - \ep^{-2} a_\t( v, \tilde{w}) - b_\t( v, \tilde{w}), \quad \forall \tilde{w} \in W.
\end{equation}
By the Cauchy-Schwarz inequality, one has
\begin{equation}\label{Nbound}
A_{\ep,\t}[\N_{\ep,\t} v] \le A_{\ep,\t}[v], \quad \forall v \in V. 
\end{equation}
One clearly has the direct sum decomposition
\[
H^1_{per}(\square) = (I+\mathcal{N}_{\ep,\theta}) V \dot{+} W,
\]
and notice that problem \eqref{p0} diagonalises with respect to this decomposition: $u_{\ep,\t} = (I+\mathcal{N}_{\ep,\theta}) v_{\ep,\t} + w_{\ep,\t}$ where $v_{\ep,\t} \in V$ solves
\begin{equation}\label{p1}
  \ep^{-2}a_{\t}\big( (I+\mathcal{N}_{\ep,\theta}) v_{\ep,\t},(I+\mathcal{N}_{\ep,\theta})\tilde{v}\big) + b_{\t}\big((I+\mathcal{N}_{\ep,\theta})v_{\ep,\t},(I+\mathcal{N}_{\ep,\theta})\tilde{v} \big) = \big(f,(I+\mathcal{N}_{\ep,\theta})\tilde{v} \big), \quad \forall \tilde{v} \in V,  
\end{equation}
and
 $w_{\ep,\t} \in W$ solves
\begin{equation}\label{Wprob.1}
  \ep^{-2}a_{\t}\big(  w_{\ep,\t},\tilde{w}\big) + b_{\t}\big(w_{\ep,\t},\tilde{w} \big) = (f,\tilde{w} ), \quad \forall \tilde{w} \in W.  
\end{equation}
Our first main result demonstrates that $w_{\ep,\t}$ is small:
\begin{lemma}
\begin{align*}
\ep^{-2}a_{\t}[  w_{\ep,\t}] + b_{\t}[w_{\ep,\t}] & \lesssim \ep^2 \| f\|_{L^2(\square)}^2, \\
\| w_{\ep,\t} \|_{L^2(\square)} & \lesssim \ep^2 \| f\|_{L^2(\square)}.
\end{align*}
\end{lemma}
\begin{proof}
Let $\tilde{w} = w_{\ep,\t}$ in \eqref{Wprob.1}. Then, by Cauchy-Schwarz, one has
\[
\ep^{-2}a_{\t}[  w_{\ep,\t}] + b_{\t}[w_{\ep,\t}] \le  \|f\|_{L^2(\square)} \|w_{\ep,\t}\|_{L^2(\square)}.
\]
Now, by Proposition \ref{prop.Wuniformcoercive}, one has
\[
\|w_{\ep,\t}\|_{L^2(\square)}^2 \lesssim a_\t[w_{\ep,\t}],
\]
and the desired inequalities follow.
\end{proof}
We have just shown that the solutions $u_{\ep,\t}$ and $v_{\ep,\t}$ to problems \eqref{p0} and \eqref{p1}, respectively, are close in the following sense:
\begin{equation*}
 \ep^{-2}a_{\t}[  u_{\ep,\t}- (I+\mathcal{N}_{\ep,\theta}) v_{\ep,\t}] + b_{\t}[u_{\ep,\t}- (I+\mathcal{N}_{\ep,\theta}) v_{\ep,\t}] \lesssim \ep^2 \|f\|_{L^2(\square)}^2,
\end{equation*}
and 
\begin{equation}
 \|  u_{\ep,\t}- (I+\mathcal{N}_{\ep,\theta}) v_{\ep,\t}\|_{L^2(\square)} \lesssim \ep^2 \|f\|_{L^2(\square)}.\label{IH1}
\end{equation}

In what follows, we will simplify further the approximate problem \eqref{p1}. First, we shall show that the term $\mathcal{N}_{\ep,\theta}\tilde{v}$ in the right-hand-side of \eqref{p1} can be neglected if we are only concerned with $O(\ep)$-approximations in $L^2$. 
To this end, let $v^{(1)}_{\ep,\t} \in V$ solve
\begin{equation}\label{p2}
  \ep^{-2}a_{\t}\big( (I+\mathcal{N}_{\ep,\theta}) v^{(1)}_{\ep,\t},(I+\mathcal{N}_{\ep,\theta})\tilde{v}\big) + b_{\t}\big((I+\mathcal{N}_{\ep,\theta})v^{(1)}_{\ep,\t},(I+\mathcal{N}_{\ep,\theta})\tilde{v} \big) = (f,\tilde{v} ), \quad \forall \tilde{v} \in V.
\end{equation}
\begin{lemma}
\begin{align}
    A_{\ep,\t}[ (I+\mathcal{N}_{\ep,\theta}) ( v_{\ep,\t}- v^{(1)}_{\ep,\t})] 
& \lesssim \ep^2 \|f\|_{L^2(\square)}^2; \nonumber \\
 \| (I+\mathcal{N}_{\ep,\theta}) ( v_{\ep,\t}- v^{(1)}_{\ep,\t})\|_{L^2(\square)}
& \lesssim \ep \|f\|_{L^2(\square)}. \label{IH2}
\end{align}
\end{lemma}
\begin{proof}
By problems \eqref{p1} and \eqref{p2}, the function $r = v_{\ep,\t}- v^{(1)}_{\ep,\t} \in V$ solves
\[
A_{\ep,\t}\big((I+\mathcal{N}_{\ep,\theta}) r,(I+\mathcal{N}_{\ep,\theta}) \tilde{v} \big) = (f, \mathcal{N}_{\ep,\t} \tilde{v}), \quad \forall \tilde{v} \in V.
\]
Thus
\[
A_{\ep,\t}[(I+\mathcal{N}_{\ep,\theta}) r] \le \|f\|_{L^2(\square)}\|\N_{\ep,\t}r\|_{L^2(\square)} \lesssim \|f\|_{L^2(\square)}\sqrt{a_\t[\N_{\ep,\t}r]} \lesssim \ep\|f\|_{L^2(\square)}\sqrt{A_{\ep,\t}[r]}, 
\]
where the last two inequalities used Proposition \ref{prop.Wuniformcoercive} and \eqref{Nbound} respectively.
Finally, by  \eqref{cor1}, one has
\[
A_{\ep,\t}[r] \lesssim A_{\ep,\t}[(I+\N_{\ep,\t})r]
\]
and the desired result follows.
\end{proof}

\section{Asymptotic analysis of correctors}
Our next task is to analyse the asymptotic behaviour of $\N_{\ep,\t}$. We will show that, up to some boundary layer, $\N_{\ep,\t}$ is asymptotically described by the usual correctors from homogenisation theory. 

We denote  by $N'$ the harmonically extended homogenisation corrector for two-dimensional perforated domains. Namely, 
the $\square'$-periodic function $N'(y')=(N_1(y'),N_2(y'))$ that satisfies

\begin{equation}\label{corp12}
\int_{\square' \backslash B} N_j = 0, \quad \text{and}\quad    \int_{\square' \backslash B} \nabla' N_j \cdot \nabla' \phi = -\int_{\square' \backslash B} e_j \cdot \nabla' \phi, \qquad \forall \phi \in H^1_{per}(\square' \backslash B),
\end{equation}
with $-\Delta' N_j = 0$ in $B$, for $j=1,2$. 

Additionally, we denote by $N_3$  the standard one-dimensional homogenisation corrector. Namely,  $N_3 \in H^1_{per}[-\tfrac{1}{2},\tfrac{1}{2})$ that  satisfies $\int_{-\tfrac{1}{2}}^{\tfrac{1}{2}} N_3(y_3) \, \mathrm{d}y_3 =0$ and solves
\begin{equation}\label{irregN3}
\int_{-\frac{1}{2}}^{\frac{1}{2}}a(y_3) \partial_3  N_3 \overline{\partial_3 \phi} \, \mathrm{d}y_3 =  - \int_{-\frac{1}{2}}^{\frac{1}{2}}a(y_3) \overline{\partial_3  \phi} \, \mathrm{d}y_3,
\qquad \text{for all $\phi \in H^1_{per}[-\tfrac{1}{2},\tfrac{1}{2})$.}    
\end{equation}
To proceed,  we write
\[
V= \CC \dot{+} H_0^1(B),
\]
where we understand $H^1_0(B)$ as a subset of $H^1_{per}(\square)$ by the obvious identification.
We shall consider the asymptotic behaviour of $\mathcal{N}_{\ep,\t}$ on $\CC$ and $H^1_0(B)$ separately.

In what follows, it will be useful to note that 
\begin{equation}\label{czgoodangle}
|c|^2 + \| z\|_{L^2(F)}^2 \le(1- |B|^{1/2})^{-1}\| c+ z\|_{L^2(\square)}^2, \qquad \forall c \in \CC, z \in H^1_0(B).
\end{equation}

\subsection{Asymptotics of $\N_{\ep,\t}$ on $H^1_0(B)$}\label{sect.NonZ}
We will show, to leading-order, that $\mathcal{N}_{\ep,\t}$ on $H^1_0(B)$ is  multiplication by a regular function.

Now consider the mapping: 
\begin{equation}\label{regN}
z \mapsto ( \mathbf{N}_{\t_3} z)(y): = \i \t_3 N_3(y_3) z(y'), \qquad z \in H^1_0(B).
\end{equation}
Clearly, for each $z$ and $\t_3$, the function $\mathbf{N}_{\t_3} z$ belongs to $W$, cf. \eqref{W}.
Furthermore, since $N_3 \in  H^1(-\frac{1}{2},\frac{1}{2})$ with $\int_{-1/2}^{1/2} N_3=0$, and $N_3$ and $z$ depending on different variables, one has 
\begin{equation}\label{Nregbds}
\|\mathbf{N}_{\t_3}z\|_{L^2(F)} + \|\partial_3 \mathbf{N}_{\t_3}z\|_{L^2(F)} \lesssim |\t_3| \|z\|_{L^2(F)}, \quad \text{and} \quad b_\t[\mathbf{N}_{\t_3}z] \lesssim|\t_3|^2 b_\t[z] \lesssim |\t_3|^2 \| \nabla ' z\|_{L^2(B)}^2,
\end{equation}
and
\begin{equation}\label{Nregbds2}
b_\t(\mathbf{N}_{\t_3}z, \tilde{c}+\tilde{z}) =0, \quad \forall \tilde{c} \in \CC, \forall z,\tilde{z} \in H^1_0(B).
\end{equation}
The following result holds:
\begin{proposition}\label{prop.z1ass}
      \begin{equation*}
      A_{\ep,\t}[\N_{\ep,\t} z - \mathbf{N}_{\t_3} z]  \lesssim \ep^2 \Big( (1 + \ep^{-2}|\t_3|^2)^2  \| z\|_{L^2(B)}^2 +  (1+\ep^{-1}|\t_3|)^2 \| \nabla' z\|_{L^2(B)}^2 \Big). 
    \end{equation*}
\end{proposition}
\begin{proof}
     Fix $w \in W$ and $z \in H^1_0(B)$.  By \eqref{firstcorrector}, $\N_{\ep,\t}z$ satisfies
     \[
     a_\t\big(\N_{\ep,\t} z, w\big) +\ep^2 b_\t\big(\N_{\ep,\t}  z,w\big ) = -\i \t_3 \int_F a(y_3) z \cdot \overline{\partial_3 w} -\i \t_3 \int_F a(y_3) z \cdot \overline{ \i \t_3 w} -\ep^2 b_\t(  z,w).
     \]
Additionally, by    the definition of $\mathbf{N}_{\t_3} z$ (cf. \eqref{regN}), we readily compute  
     \begin{flalign*}
     	a_\t(\mathbf{N}_{\t_3} z, w) 
     	& = -\i \t_3 \int_{F}  a(y_3)  z \cdot \overline{ \partial_3 w}+ \int_{F} a(y_3)  \i \t_3 \mathbf{N}_{\t_3} z \cdot \overline{\partial_3 w} +  \int_{F} a(y_3) (\partial_3 + \i \t_3) \mathbf{N}_{\t_3} z \cdot \overline{  \i \t_3 w}.
     \end{flalign*}
     So
\[
 a_\t(\N_{\ep,\t} z- \mathbf{N}_{\t_3} z, w) +\ep^2 b_\t(\N_{\ep,\t}z - \mathbf{N}_{\t_3} z,w )  
 =  R  -\ep^2 b_\t( z,w) -\ep^2 b_\t( \mathbf{N}_{\t_3} z,w),
\]
where
\[
R: =  -|\t_3|^2 \int_{F} a(y_3) z  \overline{w} -  \int_{F} a(y_3)  \i \t_3 \mathbf{N}_{\t_3} z \cdot \overline{\partial_3 w} -  \int_{F} a(y_3) (\partial_3 + \i \t_3) \mathbf{N}_{\t_3} z \cdot \overline{  \i \t_3 w}.
\]
Now, one has (cf. \eqref{Nregbds}, $|\t_3| \le \pi$ and Proposition \ref{prop.Wuniformcoercive})
 \begin{flalign*}
     \left|R \right| &  \lesssim |\t_3|^2 \|z\|_{L^2(F)}\|w\|_{L^2(F)}+ |\t_3| \|\mathbf{N}_{\t_3}z\|_{L^2(F)}\|\partial_3 w\|_{L^2(F)} + |\t_3| \|(\partial_3+\i \t_3) \mathbf{N}_{\t_3}z\|_{L^2(F)} \|w\|_{L^2(F)}\\
     & \lesssim |\t_3|^2 \|z\|_{L^2(B)}\sqrt{a_\t[w]}.
 \end{flalign*}
The fact that  $\t \in [-\pi,\pi)^3$, along with Poincar\'{e}'s inequality for $B$, gives
\begin{equation}\label{bzequivnorm}
b_\t[\tilde{z}] \lesssim \| \nabla ' \tilde{z} \|_{L^2(B)}^2, \quad \forall \tilde{z} \in H^1_0(B).
\end{equation}
This inequality, along with  Proposition \ref{prop.bcrossterm}, implies
\[
|b_\t(z, w)| \lesssim \| \nabla' z\|_{L^2(B)}\sqrt{a_\t[w]}.
\]
Finally, by Cauchy-Schwarz inequality, \eqref{bzequivnorm} and \eqref{Nregbds},  one has
\begin{flalign*}
    |b_\t(\mathbf{N}_{\t_3} z, w)|
     \lesssim |\t_3|  \| \nabla' z\|_{L^2(B)}\sqrt{b_\t[w]}.
\end{flalign*}
Therefore,
\begin{flalign*}
&\big|  a_\t(\N_{\ep,\t}z -\mathbf{N}_{\t_3} z, w) +\ep^2 b_\t(\N_{\ep,\t}z - \mathbf{N}_{\t_3} z,w ) \big| \\
& \lesssim   |\t_3|^2   \|z\|_{L^2(B)}\sqrt{a_\t[w]} + \ep^2  \|\nabla' z\|_{L^2(B)}\sqrt{a_\t[w]} + \ep^2 |\t_3|  \|\nabla' z\|_{L^2(B)}\sqrt{b_\t[w]} \\
& =  \big( |\t_3|^2   \|z\|_{L^2(B)}+ \ep^2  \|\nabla' z\|_{L^2(B)} \big) \sqrt{a_\t[w]} + \ep |\t_3|  \|\nabla' z\|_{L^2(B)}\sqrt{\ep^2 b_\t[w]} \\
&\lesssim \sqrt{|\t_3|^4 \|z\|_{L^2(B)}^2 +  (\ep^4  + \ep^2|\t_3|^2) \|\nabla' z\|_{L^2(B)}^2  }\sqrt{a_\t[w]+\ep^2 b_\t[w]} \\
& \lesssim \ep^2 \sqrt{(\ep^{-2}|\t_3|^2+1)^2 \|z\|_{L^2(B)}^2 +  (1  + \ep^{-1}|\t_3|)^2 \|\nabla' z\|_{L^2(B)}^2  }\sqrt{a_\t[w]+\ep^2 b_\t[w]}. 
\end{flalign*}
This implies the desired result.
\end{proof}

\subsection{Asymptotics of $\N_{\ep,\t}$ on $\CC$} 
Consider the $\square$-periodic function
\begin{equation*}
N_\t(y)  : = \i \t_1 N_1(y')+\i \t_2 N_2(y') + \i \t_3 \chi_0(y') N_3(y_3)  .
\end{equation*}
Notice that $\int_{\square \backslash F} N_\t  = 0$,
\[
- \Delta' \langle N_\t \rangle_3 =0, \quad \text{in $B$}, 
\]
and, by \eqref{corp12} and \eqref{irregN3}
\begin{flalign}\nonumber
    \int_{\square\backslash F} \nabla N_\t \cdot \overline{\nabla\tilde{w}} +  \int_{F}a(y_3) \partial_3  N_\t \cdot \overline{\partial_3\tilde{w}} &= -\int_{\square \backslash F} \i \t' \cdot \overline{\nabla' \tilde{w}} -     \int_F a(y_3)  \i \theta_3  \overline{\partial_3 \tilde{w}} \\
&    = -\int_{\square \backslash F} \i \t \cdot \overline{\nabla \tilde{w}} -     \int_F a(y_3)  \i \theta_3   \overline{\partial_3 \tilde{w}}, \qquad \forall \tilde{w} \in W.\label{15.10.26e1}
\end{flalign}
Furthermore, it is clear that 
\begin{equation}\label{Ntapriori}
 \|N_\t\|_{H^1(F)}+ \|N_\t\|_{H^1(\square \backslash F)}\lesssim |\t|.
\end{equation}
However,  $N_\t $ does not belong to $H^1_{per}(\square)$ unless $a$ is constant (or $\t_3=0$). We therefore introduce a `boundary layer' correction that accounts for the jump in $\chi_0 N_3$ at the boundary $\partial B \times (-\frac{1}{2},\frac{1}{2})$.
To this end, consider $\eta_\ep$ the solution to 
\begin{equation}\label{mprob1} 
-\ep^{2}\Delta' \eta_\ep - \partial_3 a(y_3) \partial_3 \eta_\ep  = 0 \quad \text{in $F$}
\end{equation}
that satisfies  $\eta_\ep(y + e_3) = \eta_\ep(y)$, and 
\begin{equation}\label{mprob2}
\eta_\ep(y) = - N_3(y_3), \quad \text{for  $y \in \partial B \times (-\tfrac{1}{2},\tfrac{1}{2})$}.    
\end{equation}
Note that  $\langle \eta_\ep \rangle_3=0$. Now, consider the $H^1_{per}(\square)$ function
\begin{equation*}
\mathbf{N}_{\ep,\t}(y) : = N_\t(y) + \i \t_3 \chi_0(y') \eta_\ep(y)  =   \i \t' \cdot N'(y') + \i \t_3\chi_0(y') \Big(N_3(y_3)+\eta_\ep(y)\Big).
\end{equation*}
Notice that $\mathbf{N}_{\ep,\t}$ belongs to $W$, cf. \eqref{W}.
The following holds:
\begin{proposition}\label{prop.truncZ0}

 \begin{equation*}
A_{\ep,\t}[\N_{\ep,\t} c - c\mathbf{N}_{\ep,\t} ] \lesssim  \ep^2 (1+\ep^{-2}|\t|^2)^2 |c|^2.
\end{equation*}
\end{proposition}
To prove Proposition \ref{prop.truncZ0} we shall use the following result, whose proof is deferred until the end of this section:
\begin{proposition}[Boundary layer estimates]\label{prop.bdry}
\begin{align}
 \label{bdyest3}
 \| \nabla' \eta_\ep \|_{L^2(F)} &\lesssim \ep^{-1/2};\\
\label{bdryest1}
\| \eta_\ep \|_{L^2(F)}+\| \partial_3 \eta_\ep \|_{L^2(F)}  &\lesssim \sqrt{\ep};\\ \label{mL1est} 
      \left| \int_F a(y_3)  \eta_\ep \right| & \lesssim \ep; \\
      \left| \int_F a(y_3)  \partial_3 \eta_\ep \right| & \lesssim \ep. \label{dmL1est} 
 \end{align}    
 Furthermore, $\partial_n \eta_\ep \in L^2\big(\partial B \times (-\frac{1}{2},\frac{1}{2})\big)$ and 
 \begin{equation}
       \label{bdryest2}
\| \partial_n \eta_\ep \|_{L^2(\partial B \times (-\frac{1}{2},\frac{1}{2}))} \lesssim \ep^{-1}.
 \end{equation}
\end{proposition}

\begin{proof}[Proof of Proposition \ref{prop.truncZ0}] Fix $c \in \CC$ and  $w \in W$,  and  set 
\[
r :=  \N_{\ep,\t} c- c \mathbf{N}_{\ep,\t} = \N_{\ep,\t} c - cN_{\t} - c\i \t_3 \chi_0 \eta_\ep .
\]
Note that $a_\t$ is well-defined on $\{ u \in L^2(\square)\, : \, \partial_3 u \in L^2(\square), \nabla' u\in L^2(\square\backslash F)\}$ and that  $b_\t$ is well-defined on $\{ u \in L^2(\square)\, : \,  \nabla'u\in L^2(F)\}$.

Now, 
\begin{flalign}
\label{18.9.26e3}    &a_\t(r,w) + \ep^2 b_\t(r,w)  \\
 \nonumber   &=  a_\t(\N_{\ep,\t} c,w) + \ep^2 b_\t(\N_{\ep,\t} c,w) - c\big(  a_\t(N_{\t},w) + \ep^2 b_\t(N_{\t},w) \big) - c\big(  a_\t(\i \t_3 \chi_0 \eta_\ep ,w) + \ep^2 b_\t(\i \t_3 \chi_0 \eta_\ep ,w) \big)\\
 \nonumber   & = -\big(  a_\t(c,w) + \ep^2 b_\t(c,w) \big)-c\big(  a_\t(N_{\t},w) + \ep^2 b_\t(N_{\t},w) \big) - c\big(  a_\t(\i \t_3 \chi_0 \eta_\ep,w) + \ep^2 b_\t(\i \t_3 \chi_0 \eta_\ep ,w) \big) \\
 \nonumber   & =-c\Big[ \underbrace{\big(  a_\t(1,w) +  a_\t(N_{\t},w)  \big)}_{=:I_1} +  \underbrace{ \big(  a_\t(\i \t_3 \chi_0 \eta_\ep ,w) + \ep^2 b_\t(\i \t_3 \chi_0 \eta_\ep ,w) \big)}_{=:I_2} +   \underbrace{\ep^2 \big( b_\t(1,w)+ b_\t( N_\t ,w) \big)}_{=:I_3} \Big] 
\end{flalign}

\textbullet\,  \textit{Estimating $I_1$.}  From the definition of $a_\t$ and \eqref{15.10.26e1}, it follows that
\begin{flalign*}
    a_\t(N_\t , w)  &= \int_{\square \backslash F} (\nabla + \i \t) N_\t  \cdot \overline{(\nabla + \i \t) w} + \int_{F} a(y_3) (\partial_3 + \i \t_3) N_\t  \cdot \overline{(\partial_3 + \i \t_3) w} \\
    & = \int_{\square \backslash F} \nabla  N_\t  \cdot \overline{\nabla  w} + \int_{F} a(y_3) \partial_3 N_\t  \cdot \overline{\partial_3  w} + R \\
    & = -\i  \Big( \int_{\square \backslash F} \t \cdot \overline{\nabla w} + \t_3 \int_F a(y_3) \overline{\partial_3 w} \Big) +R,
\end{flalign*}
where
\[
R: =  \int_{\square \backslash F} \i \t N_\t  \cdot \overline{\nabla w }+\int_{\square \backslash F} (\nabla + \i \t) N_\t  \cdot \overline{ \i \t w} + \int_{F} a(y_3)  \i \t_3 N_\t  \cdot \overline{\partial_3  w}  + \int_{F} a(y_3) (\partial_3 + \i \t_3) N_\t  \cdot \overline{ \i \t_3 w}.
\]
Therefore
\begin{flalign*}
   I_1 & =  a_\t(1,w) + a_\t(N_\t ,w) \\
   & = \int_{\square \backslash F} \i \t  \cdot \overline{(\nabla+ \i \t)w} + \int_F a(y_3) \i \t_3  \overline{(\partial_3 + \i \t_3) w} -  \int_{\square \backslash F} \i \t \cdot \overline{\nabla w} -  \int_F a(y_3) \i \t_3  \overline{\partial_3 w} + R \\
    & = \int_{\square \backslash F} \i \t  \cdot \overline{\i \t w} + \int_F a(y_3) \i \t_3  \overline{\i \t_3 w}  + R.
\end{flalign*}
Now,  Proposition \ref{prop.Wuniformcoercive} implies
\[
\left| \int_{\square \backslash F} \i \t  \cdot \overline{\i \t w} + \int_F a(y_3) \i \t_3  \overline{\i \t_3 w}  \right| \lesssim |\t|^2  \| w\|_{L^2(\square\backslash F)} +  |\t_3|^2  \| w\|_{L^2(F)}\lesssim |\t|^2  \sqrt{a_\t[w]}.
\]
Furthermore, \eqref{Ntapriori} and Proposition \ref{prop.Wuniformcoercive} imply that 
\begin{flalign*}
    |R|
    & \lesssim |\t|^2 \sqrt{a_\t[w]}.
\end{flalign*}
Whence
\begin{equation}\label{15.9.26e1}
|I_1|  \lesssim   |\t|^2 \sqrt{a_\t[w]} \le   |\t|^2 \sqrt{a_\t[w]+\ep^2b_\t[w]}.
\end{equation}
\textbullet\, \textit{Estimating $I_2$.} By \eqref{mprob1}-\eqref{mprob2}

\begin{flalign*}
 I_2&=   a_\t( \i \t_3 \chi_0 \eta_\ep ,w) + \ep^2  b_\t( \i \t_3 \chi_0 \eta_\ep ,w) \\
 & =  \i \t_3 \left(  \int_F a(y_3) \big( \partial_3 + \i \t_3 \big) \eta_\ep \cdot \overline{\big( \partial_3 + \i \t_3 \big) w}  +  \ep^2 \int_F \big( \nabla' + \i \t'  \big) \eta_\ep \cdot \overline{\big( \nabla' + \i \t'  \big) w} + \ep^2\int_F \eta_\ep \overline{w} \right)\\
 & =  \i \t_3 \left(  \int_F a(y_3)  \partial_3  \eta_\ep \cdot \overline{\partial_3  w}  +  \ep^2 \int_F \nabla'  \eta_\ep \cdot \overline{ \nabla'  w}  \right) + \i\t_3S \\
 & =   \i \t_3 \left(  \int_F \underbrace{( - \partial_3 a(y_3)  \partial_3  \eta_\ep  - \ep^2 \Delta' \eta_\ep)}_{=0} \cdot \overline{  w}  +  \ep^2 \int_{\partial B \times(-\frac{1}{2},\frac{1}{2})} \nabla'   \eta_\ep \cdot n'  \overline{  w}  \right) + \i\t_3S \\
 & =  \i \t_3   \ep^2 \int_{\partial B \times(-\frac{1}{2},\frac{1}{2})} \nabla'   \eta_\ep \cdot n'  \overline{  w}  + \i\t_3S,
\end{flalign*}
for
\begin{multline*}
    S:=  \int_{F}a(y_3) \i \t_3 \eta_\ep \cdot \overline{\partial_3 w } + \int_{F} a(y_3) (\partial_3 + \i \t_3)\eta_\ep\cdot \overline{ \i \t_3 w}  \\ + \ep^2 \int_{F} \i \t'  \eta_\ep \cdot \overline{ \nabla' w} 
    + \ep^2\int_{F} (\nabla' + \i \t') \eta_\ep \cdot \overline{ \i \t' w} +\ep^2 \int_F \eta_\ep \overline{w}.
\end{multline*}
Now, by \eqref{bdryest2}, the standard trace theorem and Proposition \ref{prop.Wuniformcoercive}, 
\begin{flalign*}
    \left| \int_{\partial B \times(-\frac{1}{2},\frac{1}{2})} \nabla'   \eta_\ep \cdot n'  \overline{  w}\right| & \le \| \partial_n \eta_\ep \|_{L^2(\partial B \times(-\frac{1}{2},\frac{1}{2}))} \| w \|_{L^2(\partial B \times(-\frac{1}{2},\frac{1}{2}))} \\
    &\lesssim \ep^{-1} \| w \|_{H^1(\square \backslash F)} \\
    &\lesssim \ep^{-1} \sqrt{a_\t[w]}.
\end{flalign*}
Also, one has 
\begin{flalign*}
    |S| & \lesssim |\t_3|   \| \eta_\ep \|_{L^2(F)} \|\partial_3 w\|_{L^2(F)} + |\t_3| (\| \eta_\ep \|_{L^2(F)}+\| \partial_3 \eta_\ep \|_{L^2(F)}) \| w \|_{L^2(F)} \\
    & +\ep^2 |\t'| \| \eta_\ep \|_{L^2(F)} \|\nabla ' w\|_{L^2(F)} + \ep^2  \Big(|\t'|\| \nabla' \eta_\ep \|_{L^2(F)}+(1+|\t'|)\| \eta_\ep \|_{L^2(F)}\Big) \|w\|_{L^2(F)}.
\end{flalign*}
So, utilising \eqref{bdyest3}, \eqref{bdryest1}, Proposition \ref{prop.Wuniformcoercive} and \eqref{btequiv} gives
\begin{flalign*}
    |S| & \lesssim |\t_3 |  \sqrt{\ep} \sqrt{a_\t[w]} + \ep^{5/2} |\t'|\|\nabla' w\|_{L^2(F)} + \ep^{3/2} |\t'| \|w\|_{L^2(F)} + \ep^{5/2} \|w\|_{L^2(F)}\\
    & \lesssim (|\t_3| \sqrt{\ep} +\ep^{3/2} |\t'| + \ep^{5/2}) \sqrt{a_\t[w]} + \ep^{5/2} |\t'| \sqrt{b_\t[w]} \\
    & \lesssim \sqrt{\ep}(|\t|+\ep)\sqrt{a_\t[w]+\ep^2b_\t[w]}.
\end{flalign*}
Therefore
\begin{flalign}
|I_2|  & \lesssim |\t_3|\ep^2 \ep^{-1} \sqrt{a_\t[w]} +|\t_3| \sqrt{\ep}(|\t|+\ep)\sqrt{a_\t[w]+\ep^2b_\t[w]} \nonumber \\ 
& \lesssim   |\t_3| (\ep + \sqrt{\ep}|\t|) \sqrt{a_\t[w]+\ep^2b_\t[w]}  \nonumber
 \\& \lesssim (|\t|^2+\ep^2)   \sqrt{a_\t[w]+\ep^2b_\t[w]}.\label{18.9.26e1}
\end{flalign}
\textbullet\, \textit{Estimating $I_3$.} 
Note that, by Proposition \ref{prop.bcrossterm}, one has  
\[
|b_\t(1,w)| \lesssim \sqrt{a_\t[w]}.
\]
Now Cauchy-Schwarz inequality, along with \eqref{Ntapriori} gives
\[
|b_\t(N_\t ,w)| \le \sqrt{b_\t[N_\t ]} \sqrt{b_\t[w]}\lesssim |\t|  \sqrt{b_\t[w]}.
\]
So
\begin{equation}\label{18.9.26e2}
|I_3| \lesssim  \big( \ep^2\sqrt{a_\t[w]} + |\t|\ep  \sqrt{\ep^2 b_\t[w]} \big) \lesssim (\ep^2+\ep|\t|)\sqrt{a_\t[w]+\ep^2b_\t[w]}. 
\end{equation}
Then, \eqref{18.9.26e3}, 
 \eqref{15.9.26e1}, \eqref{18.9.26e1} and \eqref{18.9.26e2} give the desired result, upon setting $w=r$.
\end{proof}

\begin{proof}[Proof of Proposition \ref{prop.bdry}]
 Denote $G := \partial B \times(-1/2,1/2)$ with outer unit normal $n(y)=\big(n'(y'),0\big)$  where $n'$ is normal to $\partial B$.    The weak form of \eqref{mprob1}-\eqref{mprob2} is 
    \begin{equation}\label{me2}
\ep^2 \int_F \nabla' \eta_\ep \cdot \nabla' \phi +  \int_F a(y_3) \partial_3 \eta_\ep \partial_3 \phi = \ep^2\int_{G} \partial_{n} \eta_\ep \phi , \quad \forall \phi \in H^1(F)\text{ that are periodic in $y_3$}.
\end{equation}
Let $d: \overline{B} \rightarrow [0,\infty)$ denote the distance function from $\partial B$:
	\[
	d(x') := \min_{y' \in \partial B} |x'-y'|.
	\]
Note that  $d(x') =0$ for $x' \in \partial B$,
 $| \nabla' d(x') | =1 \quad \text{a.e. $x' \in B$},$
and that
\begin{equation}\label{distbdry}
\nabla' d(x') = - n'(x'), \quad x' \in \partial B.    
\end{equation}
It is clear that
\begin{equation*}
\int_F e^{-\frac{d(x')}{\delta}} \, \mathrm{d}x \lesssim \delta, \quad \forall \delta>0.
\end{equation*}
From the identity
\[
\nabla' \eta_\ep\cdot \nabla' (e^{ 2d/\delta} \eta_\ep) = |\nabla' (e^{d/\delta}\eta_\ep)|^2 - \delta^{-2} \big(e^{d/\delta}\eta_\ep\big)^2,
\]
and  taking the test function $\phi = e^{ 2d/ \delta} \eta_\ep$ in \eqref{me2}, recalling \eqref{mprob2}, gives
\begin{equation*}
\ep^2 \int_F |\nabla'(e^{d/\delta}\eta_\ep) |^2 -  \frac{\ep^2}{\delta^2}\int_F  |e^{d/\delta}\eta_\ep|^2+  \int_F a(y_3) |\partial_3 ( e^{d/\delta}\eta_\ep)|^2 = -\ep^2\int_{G} \partial_n \eta_\ep N_3
\end{equation*}
Clearly, there exists $\mu >0$ such that
\[
2\mu^2 \int_{-\frac{1}{2}}^{\frac{1}{2}} |u|^2 \le  \int_{-\frac{1}{2}}^{\frac{1}{2}} a |\partial_3 u|^2, \quad \text{$u\in 
H^1_{per}(-\tfrac{1}{2},\tfrac{1}{2})$ such that $\int_{-\frac{1}{2}}^{\frac{1}{2}} u =0$.}
\]
Thus, setting  $\delta =  \ep / \mu$ gives
\[
\ep^2 \int_F |\nabla'(e^{\mu d/\ep}\eta_\ep) |^2 +\int_F  |e^{\mu d/\ep}\eta_\ep|^2 +  \int_F a(y_3) |\partial_3 ( e^{\mu d/\ep}\eta_\ep)|^2\lesssim -\ep^2\int_{G} \partial_n \eta_\ep N_3.
\]

Setting $\phi = N_3(y_3)$ in \eqref{me2} gives
  \begin{equation*}
  \int_F a(y_3) \partial_3 \eta_\ep \partial_3 N_3 = \ep^2\int_{G} \partial_n \eta_\ep N_3.
\end{equation*}
So
\begin{flalign*}
  &  \ep^2 \int_F |\nabla'(e^{\mu d/\ep}\eta_\ep) |^2 +\int_F  |e^{\mu d/\ep}\eta_\ep|^2 +  \int_F a(y_3) |\partial_3 ( e^{\mu d/\ep}\eta_\ep)|^2 \lesssim -\int_F a(y_3) \partial_3 \eta_\ep \partial_3 N_3  \\
  & = -\int_F a(y_3) \partial_3 (e^{\mu d/\ep}\eta_\ep) \partial_3 (e^{-\mu d/\ep}N_3) 
\end{flalign*}
This implies
\begin{flalign*}
&\ep^2 \int_F |\nabla'(e^{\mu d/\ep}\eta_\ep) |^2 +\int_F  |e^{\mu d/\ep}\eta_\ep|^2 +  \int_F a(y_3) |\partial_3 ( e^{\mu d/\ep}\eta_\ep)|^2  \lesssim
  \int_F e^{-\mu 2d(x')/\ep}a(y_3)|\partial_3 N_3(y_3)|^2 \\
  &=   \int_B e^{-2\mu d(x')/\ep}\int_{-1/2}^{1/2} a(y_3)|\partial_3 N_3(y_3)|^2     \\
  & \lesssim \ep.
\end{flalign*}
Therefore, \eqref{bdryest1} and \eqref{bdyest3} readily follow.

Now
\begin{flalign*}
 \left|   \int_F a(y_3) \eta_\ep \right|& =\left|     \int_F a(y_3) e^{\mu d/\ep}\eta_\ep e^{-\mu d/\ep} \right|\\
&\le \left( \int_F  |e^{\mu d/\ep}\eta_\ep|^2\right)^{1/2} \left( \int_F  a^2(y_3) e^{-\mu 2d/\ep}\right)^{1/2} \\
& = \left( \int_F  |e^{\mu d/\ep}\eta_\ep|^2\right)^{1/2} \left( \int_{-1/2}^{1/2}  a^2 \right)^{1/2}\left( \int_B  e^{-\mu 2d/\ep}\right)^{1/2} \\
& \lesssim \ep.
\end{flalign*}
Therefore \eqref{mL1est} holds.  Similarly, 
\eqref{dmL1est} holds.

\textit{Proof of \eqref{bdryest2}.} Since $d$ is only regular in a neighbourhood of the boundary $\partial B$, let us take $\tau$ to be the cut-off function for a neighbourhood of the boundary where $d$ is regular and consider $D = \tau d$.

We shall now take $\phi = \nabla' D \cdot \nabla' \eta_\ep$ as a test function in \eqref{me2}.\footnote{ This is legitimate when $a$ is smooth, as $\eta_\ep$ is then regular near $G$. For general $a$, one first replaces $a$ by a mollification, which still satisfies $\nu<a<\nu^{-1}$, and then passes to the limit. This is possible because all the constants below depend on $a$ only through $\nu$.} Upon noting that  $\phi =- \partial_n \eta_\ep$ on $G$, cf \eqref{distbdry}, gives
\begin{equation}
\ep^2 \int_{G} |\partial_n \eta_\ep|^2= -  \ep^2 \int_F \nabla' \eta_\ep \cdot \nabla' \phi -  \int_F a(y_3)\partial_3 \eta_\ep \partial_3 \phi. \label{19.10.26e1}    
\end{equation}
To calculate the first integral on the right, we note that
\begin{flalign*}
    \nabla' \eta_\ep \cdot \nabla' \phi  =(\nabla'^2 D) \nabla' \eta_\ep \cdot \nabla' \eta_\ep  +  \frac{1}{2}  \nabla' D \cdot \nabla' \big(|\nabla' \eta_\ep|^2 \big).
\end{flalign*}
So, one has
\begin{flalign*}
\ep^2 \int_F \nabla' \eta_\ep \cdot \nabla' \phi & =  \ep^2 \int_F (\nabla'^2 D) \nabla' \eta_\ep \cdot \nabla' \eta_\ep  +   \ep^2 \int_F \frac{1}{2}  \nabla' D \cdot \nabla' \big(|\nabla' \eta_\ep|^2 \big) \\
& =  \ep^2 \int_F (\nabla'^2 D) \nabla' \eta_\ep \cdot \nabla' \eta_\ep  - \ep^2 \int_F \frac{1}{2}  (\Delta' D) |\nabla' \eta_\ep|^2 + \ep^2 \int_{G} \frac{1}{2}  \nabla' D \cdot n' |\nabla' \eta_\ep|^2 \\
& =  \ep^2 \int_F (\nabla'^2 D) \nabla' \eta_\ep \cdot \nabla' \eta_\ep  - \ep^2 \int_F \frac{1}{2} ( \Delta' D )|\nabla' \eta_\ep|^2  - \frac{1}{2} \ep^2 \int_{G}  |\nabla' \eta_\ep|^2.
\end{flalign*}
Also, for the second integral, notice that since $\partial_3 D= 0$, then
\[
  \partial_3 \eta_\ep \partial_3 \phi =  \frac{1}{2}  \nabla' D \cdot \nabla' \big(|\partial_3 \eta_\ep|^2 \big)
\]
So, one has
\begin{flalign*}
 \int_F a(y_3)\partial_3 \eta_\ep \partial_3 \phi & =  \int_F a(y_3)  \frac{1}{2}  \nabla' D \cdot \nabla' \big(|\partial_3 \eta_\ep|^2 \big) \\
& =  - \int_F a(y_3)  \frac{1}{2} ( \Delta' D )|\partial_3 \eta_\ep|^2  + \int_{G} a(y_3)  \frac{1}{2}  \nabla' D \cdot n' \big(|\partial_3 \eta_\ep|^2 \big)\\
& =   -  \frac{1}{2} \int_F a(y_3)  (\Delta' D )|\partial_3 \eta_\ep|^2  -  \frac{1}{2} \int_{G} a(y_3)  |\partial_3 \eta_\ep|^2.
\end{flalign*}
Now, since  $\eta_\ep(y',y_3) = -N_3(y_3)$ on $G$, then $\partial_3 \eta_\ep = -\partial_3 N_3$ and $\nabla' \eta_\ep = (\partial_n \eta_\ep) n'$ on $G$. Therefore, \eqref{19.10.26e1} can be rewritten as \begin{flalign*}
  & \ep^2 \int_{G} |\partial_n \eta_\ep|^2 =    - 2\ep^2 \int_F (\nabla'^2 D) \nabla' \eta_\ep \cdot \nabla' \eta_\ep  + \ep^2 \int_F  ( \Delta' D )|\nabla' \eta_\ep|^2 +   \int_F a ( \Delta' D)|\partial_3 \eta_\ep|^2  + \int_{G} a  |\partial_3 N_3|^2   \\
   & \lesssim   \ep^2 \| \nabla' \eta_\ep\|_{L^2(F)}^2  + \| \partial_3 \eta_\ep \|_{L^2(F)}^2 +1 \lesssim    1. 
\end{flalign*}
That is,  \eqref{bdryest2}  holds.
\end{proof}

\section{The limit form}\label{sec.limform}
In the later sections, we will establish the leading-order quantitative estimates for the solution to problem \eqref{p2}. This section is devoted to the introduction of the limit objects.

Let us introduce the limit form on $V=\CC\dot{+}H^1_0(B)$:
\begin{flalign*}
 &   A^h_{\ep,\t}(c+z,\tilde{c}+\tilde{z}) \\
    & := \ep^{-2} A^h \t \cdot \t  c \overline{\tilde{c}}    + \ep^{-2} a^h|\t_3|^2 (c+z,\tilde{c}+\tilde{z})_{L^2(B)}+   \big((\nabla'+\i\t')z,(\nabla'+\i\t')\tilde{z}\big)_{L^2(B)}+
   (c+z,\tilde{c}+\tilde{z}),
\end{flalign*}
for $c,\tilde{c} \in \CC$ and $z,\tilde{z} \in H^1_0(B)$. Here,
 $a^{\mathrm h}: = \left(\int_{-1/2}^{1/2} a^{-1} \right)^{-1}$ is the one-dimensional homogenised coefficient and $A^{\mathrm h}$ is the constant, symmetric, positive-definite $3\times 3$  perforated-fibre homogenised matrix given, in \cite[(2.21)]{ChSmZh}, as 
\[
A^{\mathrm h} : = |\square\backslash F| e_3\otimes e_3 +  \sum_{\alpha,\beta=1}^2 e_\alpha \otimes e_\beta\int_{\square'\backslash B} (\partial_{\beta} N_\alpha + \delta_{\alpha \beta}),
\] 
where $N_\alpha$ are the homogenisation correctors given by \eqref{corp12}.

Now, we consider the limit problem: Find $c+z$, $c \in \CC$ and $z \in H^1_0(B)$ the solution to 
\begin{equation}\label{prob.hom}
A^h_{\ep,\t}(c+z,\tilde{c}+\tilde{z}) = (f,\tilde{c} +\tilde{z}), \quad \forall \tilde{c} \in \CC, \forall \tilde{z} \in H^1_0(B).   
\end{equation}
The following a priori estimates hold.
\begin{proposition}\label{prop.homapriori}
    The solution $c+ z$ to  \eqref{prob.hom} satisfies
 \begin{align}
    \big(1+ \ep^{-2} |\t_3|^2\big) \| z \|_{L^2(B)} & \lesssim \| f\|_{L^2(\square)}; \label{8.10.26e1}\\
     \big(1+ \ep^{-1} |\t_3|\big) \|  \nabla' z \|_{L^2(B)} & \lesssim \| f\|_{L^2(\square)}; \label{8.10.26e2}\\
       \big(1+ \ep^{-2} |\t|^2\big) |c| & \lesssim \| f\|_{L^2(\square)}. \label{8.10.26e3}
 \end{align} 
  \end{proposition}
\begin{proof}

\textbullet\,     Let $\tilde{c}= c$ and $\tilde{z}=  z$ in \eqref{prob.hom}, then
    \begin{equation}\label{8.10.26e1.1}
A^h_{\ep,\t}[c+z] \le \| f\|_{L^2(\square)} \| c+z\|_{L^2(\square)},
    \end{equation}
and  it clearly  follows that 
\begin{equation}\label{sumbd9.10.26}
  \big(1+ \ep^{-2} |\t_3|^2\big) \| c+ z \|_{L^2(\square)}  \lesssim \| f\|_{L^2(\square)},  
\end{equation}
which, along with \eqref{czgoodangle}, gives \eqref{8.10.26e1}.

\textbullet\, Now, from \eqref{8.10.26e1.1} and \eqref{sumbd9.10.26}, one has 
\begin{flalign*}
 \| \nabla' z + \i \t' z\|_{L^2(B)}^2 & \le \| f\|_{L^2(\square)} \| c+z \|_{L^2(\square)} \\
    & \lesssim (1+\ep^{-2} |\t_3|^2)^{-1} \| f\|_{L^2(\square)}^2.
\end{flalign*}
It remains to note that, since $e^{\i\t'\cdot y'}z =0$ on $\partial B$, by the Poincar\'{e} inequality, with constant $C_p$, and the fact that $|\t'| \le \sqrt{2}\pi$, one clearly has 
\[
\| \nabla' z\|_{L^2(F)} \le  \| \nabla' z + \i \t' z\|_{L^2(F)}+  |\t'| \|e^{\i\t'\cdot y'}z\|_{L^2(F)} \le (1+\sqrt{2}\pi C_P) \| \nabla' z + \i \t' z\|_{L^2(F)},
\]
and \eqref{8.10.26e2} follows.

\textbullet\, 
Taking $\tilde{c} = c$ and $\tilde{z}=0$, in \eqref{prob.hom}, and utilising \eqref{8.10.26e1}, gives
\begin{flalign*}
    A^h_{\ep,\t}[c] &= (f,c) -A^h_{\ep,\t}(z,c) \\
&= (f,c) - (\ep^{-2}a^h |\t_3|^2 + 1) (z,c)_{L^2(B)} \\
& \lesssim \| f\|_{L^2(\square)} |c| + (1+ \ep^{-2}|\t_3|^2) \| z\|_{L^2(B)} |c| \\
& \lesssim \| f\|_{L^2(\square)} |c|,
\end{flalign*}
and, since  one clearly has 
\[
(1+\ep^{-2}|\t|^2) |c|^2 \lesssim A^h_{\ep,\t}[c],
\]then \eqref{8.10.26e3} follows.
\end{proof}
\section{Reduction to the limit form}
Here, we shall prove the following result.
\begin{theorem}\label{t.fin}
Consider $v^{(1)}_{\ep,\t}$ and $v^{(2)} : = c+z$ the solutions to \eqref{p2} and  \eqref{prob.hom}, respectively. Then
\begin{align}
A_{\ep,\t}[ (I+\mathcal{N}_{\ep,\theta}) (v^{(1)}_{\ep,\t}-v^{(2)})] &\lesssim \ep^2 \|f\|_{L^2(\square)}^2; \label{finalest1}\\
 \| (I+\mathcal{N}_{\ep,\theta})  v^{(1)}_{\ep,\t}- v^{(2)}\|_{L^2(\square)}
& \lesssim \ep \|f\|_{L^2(\square)}.\label{IH3}
\end{align}
\end{theorem}
To prove this, we shall utilise the following result, whose proof is delayed until the end of the section.
\begin{proposition}\label{prop.finalrelbd}
    \begin{multline*}
        \left| A_{\ep,\t}\big( (I+\N_{\ep,\t})(c+z),(I+\N_{\ep,\t})(\tilde{c}+\tilde{z})\big) - A^h_{\ep,\t}(c+z,\tilde{c}+\tilde{z})\right| 
        \lesssim \ep \Phi(c,z)  \sqrt{A_{\ep,\t}[(I+\N_{\ep,\t})(\tilde{c}+\tilde{z})]},
    \end{multline*}
    for
    \[
\Phi(c,z) = \Big( (1+\ep^{-2}|\t|^2)|c| +(1+\ep^{-2}|\t_3|^2)\| z\|_{L^2(B)} +(1+\ep^{-1}|\t_3|)\| \nabla' z\|_{L^2(B)}   \Big).
    \]
\end{proposition}

\begin{proof}[Proof of Theorem \ref{t.fin}]
The function $ r= v^{(1)}_{\ep,\t} - v^{(2)} = r_c+r_z$ for $r_c \in \CC$ and $r_z\in H^1_0(B)$ satisfies
\[
   A_{\ep,\t}[(I+\N_{\ep,\t})r]  = A^h_{\ep,\t}(c+z,r_c+r_z)-A_{\ep,\t}\big( (I+\N_{\ep,\t})(c+z), (I+\N_{\ep,\t})r \big) 
\]
which along with Proposition \ref{prop.finalrelbd} gives 
\begin{flalign*}
 &   A_{\ep,\t}[(I+\N_{\ep,\t})r]  
  \lesssim  \ep \Phi(c,z) \sqrt{A_{\ep,\t}[(I+\N_{\ep,\t})r]}.
\end{flalign*}
Finally, from Proposition \ref{prop.homapriori}, one has 
\[
\Phi(c,z)
\lesssim \| f\|_{L^2(\square)}.
\]
Hence, \eqref{finalest1} holds. 

Next, by \eqref{Nbound}, \eqref{eqnatv1}, \eqref{czgoodangle} and Proposition \ref{prop.homapriori}:
\begin{flalign*}
\| \N_{\ep,\t} v^{(2)} \|_{L^2(\square)}^2  &\lesssim a_{\t}[ \N_{\ep,\t} v^{(2)}]\\
&\le \ep^2A_{\ep,\t}[ \N_{\ep,\t} v^{(2)}]\\
&\le \ep^2A_{\ep,\t}[ v^{(2)}]\\
& \lesssim (\ep^2+|\t|^2)|c|^2 + (\ep^2+|\t_3|^2 )\| z\|_{L^2(B)}^2 + \ep^2 \| \nabla' z\|_{L^2(B)}^2 \\
& \lesssim \ep^2 \| f\|_{L^2(\square)}^2.
\end{flalign*}
Then \eqref{IH3} follows from the above assertion and \eqref{finalest1}.
\end{proof}
\begin{proof}[Proof of Proposition \ref{prop.finalrelbd}]
Throughout the proof, we shall use the following:
\begin{multline}\label{gatherkeystuff}
     (1+\ep^{-2}|\t|^2) |\tilde{c}|^2+ (1+\ep^{-2} |\t_3|^2) \| \tilde{z} \|_{L^2(B)}^2 + \| \nabla' \tilde{z} \|_{L^2(B)}^2 \lesssim A_{\ep,\t}[(I+\N_{\ep,\t})(\tilde{c}+\tilde{z})], \\ \quad \forall \tilde{c} \in \CC, \forall \tilde{z} \in H^1_0(B), \forall \t \in [-\pi,\pi)^3, \forall 0<\ep<1.
\end{multline}
Indeed, this follows from \eqref{czgoodangle},  Corollary \ref{corAonV}, Proposition \ref{prop.bsplit} and definition \eqref{A}.

To begin,  notice that the ranges of  $I+\N_{\ep,\t}$ and $\N_{\ep,\t}$ are orthogonal with respect to $A_{\ep,\t}$; thus
\begin{flalign*}
   &    A_{\ep,\t}\big( (I+\N_{\ep,\t})(c+z), (I+\N_{\ep,\t})(\tilde{c}+\tilde{z})\big)  =   A_{\ep,\t}\big( (I+\N_{\ep,\t})(c+z), \tilde{c}+\tilde{z}\big) \\
       & = A_{\ep,\t}\big((I+\N_{\ep,\t})c,\tilde{c}\big) +A_{\ep,\t}\big((I+\N_{\ep,\t})c,\tilde{z}\big) +A_{\ep,\t}\big((I+\N_{\ep,\t})z,\tilde{c} +\tilde{z}\big) 
\end{flalign*}
Furthermore,  $\N_{\ep,\t}$ is symmetric with respect to $A_{\ep,\t}$, and  so
\begin{flalign*}
A_{\ep,\t}\big((I+\N_{\ep,\t})c,\tilde{z}\big) 
    = A_{\ep,\t}\big(c,(I+\N_{\ep,\t})\tilde{z}\big).
\end{flalign*}
Hence,
\begin{multline*}
  A_{\ep,\t}\big( (I+\N_{\ep,\t})(c+z), (I+\N_{\ep,\t})(\tilde{c}+\tilde{z})\big) =  A_{\ep,\t}\big((I+\N_{\ep,\t})c,\tilde{c}\big) +A_{\ep,\t}\big(c,(I+\N_{\ep,\t})\tilde{z}\big) \\ +A_{\ep,\t}\big((I+\N_{\ep,\t})z,\tilde{c}+\tilde{z}\big).
\end{multline*}

\textbullet \, Here, we shall establish
\begin{multline}
       \Big| A_{\ep,\t}\big((I+\N_{\ep,\t})z,\tilde{c}+\tilde{z}\big) - \Big( \ep^{-2}  a^h |\t_3|^2 ( z,\tilde{c}+\tilde{z})+(z,\tilde{c}) + b_\t(z,\tilde{z})\Big) \Big|  \\ 
        \lesssim   \ep\Big( (1+\ep^{-2}|\t_3|^2) \| z\|_{L^2(B)}+ (1+\ep^{-1}|\t_3|) \|\nabla' z \|_{L^2(B)}  \Big) \sqrt{A_{\ep,\t}[(I+\N_{\ep,\t})(\tilde{c}+\tilde{z})]}.\label{9.10.26approx1}
\end{multline}
First, one has, cf. \eqref{cor1} and Proposition \ref{prop.z1ass}, 
\begin{flalign}
&\left| A_{\ep,\t}\big((I+\N_{\ep,\t})z,\tilde{c}+ \tilde{z}\big) - A_{\ep,\t}\big((I+\mathbf{N}_{\t_3})z,\tilde{c}+\tilde{z}\big) \right|  =\left| A_{\ep,\t}\big(\N_{\ep,\t}z - \mathbf{N}_{\t_3}z,\tilde{c}+\tilde{z}\big) \right| \nonumber \\
& \le \sqrt{A_{\ep,\t}[\N_{\ep,\t} z-\mathbf{N}_{\t_3}z]}\sqrt{A_{\ep,\t}[\tilde{c}+\tilde{z}]} \nonumber \\
& \lesssim \sqrt{A_{\ep,\t}[\N_{\ep,\t} z-\mathbf{N}_{\t_3}z]}\sqrt{A_{\ep,\t}[(I+\N_{\ep,\t})(\tilde{c}+\tilde{z})]} \nonumber \\
&\lesssim \ep\Big( (1+\ep^{-2}|\t_3|^2) \| z \|_{L^2(B)} + (1+\ep^{-1}|\t_3|) \| \nabla' z \|_{L^2(B)} \Big)\sqrt{A_{\ep,\t}[(I+\N_{\ep,\t})(\tilde{c}+\tilde{z})]} \label{24.10.26e1}
\end{flalign}
Next,
\begin{flalign*}
    a_\t\big((I+\mathbf{N}_{\t_3})z,\tilde{c}+\tilde{z}\big)  & =\int_F a \big( \i \t_3 z + \partial_3 \i \t_3 N_3 z + (\i\t_3)^2 N_3 z \big) \cdot \overline{\i \t_3 (\tilde{c} + \tilde{z})} \\
     & =\int_F a \big( \i \t_3 z + \partial_3 \i \t_3 N_3 z  \big) \cdot \overline{\i \t_3 (\tilde{c} + \tilde{z})} +  \int_F a \big( (\i\t_3)^2 N_3 z \big) \cdot \overline{\i \t_3 (\tilde{c} + \tilde{z})} \\
    & = a^h |\t_3|^2\int_B  z\overline{(\tilde{c}+\tilde{z})}  + \i \t_3 |\t_3|^2 \int_F a N_3 z \overline{(\tilde{c}+\tilde{z})}.
\end{flalign*}
So, 
\begin{flalign}
&\left| \ep^{-2}a_\t\big((I+\mathbf{N}_{\t_3})z,\tilde{c}+\tilde{z}\big)  - \ep^{-2}a^h |\t_3|^2 (z,\tilde{c}+ \tilde{z}) \right|  \lesssim \ep^{-1}|\t_3|^2 \| z\|_{L^2(B)} \big( \ep^{-1} |\t_3| \| \tilde{c}+\tilde{z}\|_{L^2(B)} \big)
.\label{12.10.26e1}
\end{flalign}
Also, cf. \eqref{Nregbds2},
\begin{flalign}
    \left| b_\t\big((I+\mathbf{N}_{\t_3})z,\tilde{c}+\tilde{z}\big)  - \Big( (z,\tilde{c}) + b_\t(z,\tilde{z}) \Big)  \right|
    & = \left| \int_F (\nabla'+\i \t')z \cdot \overline{\i\t' \tilde{c}} + b_\t\big( \mathbf{N}_{\t_3}z,\tilde{c}+\tilde{z}\big) \right| \nonumber \\
    &= \left| \int_F (\nabla'+\i \t')z \cdot \overline{\i\t' \tilde{c}}  \right| \nonumber\\
       &\lesssim |\t'|\| \nabla' z\|_{L^2(B)} |\tilde{c}|.
    \label{24.10.26e2}
\end{flalign}
Thus \eqref{9.10.26approx1} follows from \eqref{24.10.26e1}-\eqref{24.10.26e2} and \eqref{gatherkeystuff}.

\textbullet\, Here, we will show
\begin{multline} \label{diagmm}
    \left| A_{\ep,\t}\big(c,(I+\N_{\ep,\t})\tilde{z}\big) - \Big(\ep^{-2}a^h |\t_3|^2 (c,\tilde{z}) + (c,\tilde{z}) \Big)  \right| \\ \lesssim \ep(1+\ep^{-2}|\t|^2) |c| \sqrt{A_{\ep,\t}[(I+\N_{\ep,\t})(\tilde{c}+\tilde{z})]}.
\end{multline}
By Proposition \ref{prop.z1ass},
\begin{flalign}
&\left| A_{\ep,\t}\big(c,(I+\N_{\ep,\t})\tilde{z}\big) - A_{\ep,\t}\big(c,(I+\mathbf{N}_{\t_3})\tilde{z}\big)  \right|    =\left| A_{\ep,\t}\big(c,\N_{\ep,\t}\tilde{z} -\mathbf{N}_{\t_3}\tilde{z}\big)  \right| \nonumber \\
& \lesssim \ep \sqrt{A_{\ep,\t}[c]}\Big( (1+\ep^{-2}|\t_3|^2) \| \tilde{z} \|_{L^2(B)} + (1+\ep^{-1}|\t_3|) \| \nabla' \tilde{z} \|_{L^2(B)}\Big) \nonumber\\
& \lesssim \ep \sqrt{A_{\ep,\t}[c]}(1+\ep^{-1}|\t_3|)\Big( (1+\ep^{-1}|\t_3|)\| \tilde{z}\|_{L^2(B)}+ \| \nabla' \tilde{z} \|_{L^2(B)}\Big) \nonumber\\
& \lesssim  \ep (1+\ep^{-2}|\t|^2)|c|\Big( (1+\ep^{-1}|\t_3|)\| \tilde{z}\|_{L^2(B)}+ \| \nabla' \tilde{z} \|_{L^2(B)}\Big). \label{24.10.26e3}
\end{flalign}
 Now, by \eqref{12.10.26e1}, 
\begin{flalign}
    \left| \ep^{-2}a_\t\big(c,(I+\mathbf{N}_{\t_3})\tilde{z}\big) -\ep^{-2}a^h |\t_3|^2 (c,\tilde{z})\right| 
     & \lesssim  \ep^{-1}|\t_3|^2 |c|\big(\ep^{-1}|\t_3|\| \tilde{z} \|_{L^2(B)}\big)
      .\label{24.10.26e4}
\end{flalign}
Also, by $\t' \in[-\pi,\pi)^2$, \eqref{Nregbds2} and \eqref{cor1},
\begin{flalign}
    \left| b_\t\big(c,(I+\mathbf{N}_{\t_3})\tilde{z}\big) - (c,\tilde{z}) \right| & = \left| \int_F \i \t' c \cdot\overline{(\nabla'+\i\t')\tilde{z}}+  b_\t\big(c,\mathbf{N}_{\t_3}\tilde{z}\big) \right|= \left| \int_F \i \t' c \cdot\overline{(\nabla'+\i\t')\tilde{z}} \right| \nonumber \\
     & \lesssim |\t||c|  \| \nabla'\tilde{z}\|_{L^2(B)}
 .\label{24.10.26e5}
\end{flalign}
Thus \eqref{diagmm} follows from \eqref{24.10.26e3}-\eqref{24.10.26e5} and \eqref{gatherkeystuff}.


\textbullet\, Here, we demonstrate
\begin{multline}\label{12.10.26e3}
\left |A_{\ep,\t}\big((I+\N_{\ep,\t})c,\tilde{c}\big) - \ep^{-2}\Big( A^h \t \cdot \t  + |B| a^h |\t_3|^2 \Big)  c \overline{\tilde{c}}     - c \overline{\tilde{c}} \right|  \\ \lesssim    \ep(1+\ep^{-2} |\t|^2) |c|\sqrt{A_{\ep,\t}[(I+\N_{\ep,\t})(\tilde{c}+\tilde{z})]}.
\end{multline}
Indeed, by Proposition \ref{prop.truncZ0},
\begin{flalign}
 |A_{\ep,\t}\big((I+\N_{\ep,\t})c,\tilde{c}\big)      -A_{\ep,\t}\big(c+c\mathbf{N}_{\ep,\t},\tilde{c}\big)  |  & = |A_{\ep,\t}\big(\N_{\ep,\t}c- \mathbf{N}_{\ep,\t}c,\tilde{c}\big)  | \nonumber  \\
 & \le \sqrt{A_{\ep,\t}[\N_{\ep,\t}c- \mathbf{N}_{\ep,\t}c]} \sqrt{A_{\ep,\t}[\tilde{c}]}\nonumber  \\
 & \lesssim  \ep (1+\ep^{-2}|\t|^2) |c|\sqrt{A_{\ep,\t}[\tilde{c}]} \nonumber\\
 & \lesssim\ep (1+\ep^{-2}|\t|^2) |c|(1+\ep^{-1} |\t|)|\tilde{c}|.\label{24.10.26e6}
\end{flalign}
By utilising Proposition \ref{prop.bdry}, one has
\begin{flalign}
 |  \ep^{-2}a_{\t}\big(c+c\mathbf{N}_{\ep,\t},\tilde{c}\big) - \ep^{-2} a_{\t}\big(c+cN_{\t},\tilde{c}\big)| 
& = \ep^{-2}| \t_3|^2 \left|  \int_F a(y_3) (\partial_3 + \i \t_3) \eta_\ep \right|  |c||\tilde{c}| \nonumber\\
&\lesssim \ep^{-1}|\t_3|^2|c|  |\tilde{c}|
. \label{24.10.26e7}
\end{flalign}

Finally,
one has 
\begin{flalign*}
 &  a_\t\big(c+cN_{\t},\tilde{c}\big)  \\
 & = \int_{\square \backslash F}  \t  c \cdot \overline{\t \tilde{c}} + \int_{\square \backslash F}  \nabla'  ( \theta' \cdot  N') c \cdot \overline{ \t' \tilde{c}} + \int_{\square \backslash F}   \t ( \i \t'\cdot N') c \cdot \overline{ \t \tilde{c}}+ \int_F a\big(\i \t_3  +  \i \theta_3 \partial_3 N_3 + \i \t_3  N_\t  \big)c \cdot \overline{\i \t_3 \tilde{c}} \\
 & = \Big(  A^h \t \cdot \t  + |B| a^h |\t_3|^2 \Big)  c \overline{\tilde{c}} +   \int_{\square \backslash F}  \t (\i \t'\cdot N') c \cdot \overline{ \t \tilde{c}}+ \int_F a  \t_3 N_\t c  \overline{ \t_3 \tilde{c}}
\end{flalign*}
So,
\begin{flalign}
    \left|  \ep^{-2}a_\t\big(c+cN_{\t},\tilde{c}\big) - \ep^{-2}\Big( A^h \t \cdot \t  + |B| a^h |\t_3|^2 \Big)  c \overline{\tilde{c}} \right|
     & \lesssim \ep^{-1} |\t|^2 |c|\big(\ep^{-1} |\t| | \tilde{c}| \big)
    . \label{24.10.26e8}
\end{flalign}
Also,  since $\langle \eta_\ep \rangle_3=0$ one has
\begin{flalign}
    b_\t\big(c+c\mathbf{N}_{\ep,\t},\tilde{c}\big) -  b_\t\big(c+cN_{\t},\tilde{c}\big)
=     \i \t_3 \int_F (\nabla' +  \i \t') \eta_\ep  \cdot \overline{\i \t'}  c\overline{\tilde{c}} +\i \t_3 \int_F \eta_\ep   c\overline{\tilde{c}} = 0.
\label{24.10.26e9}
\end{flalign}
Finally, by $\t'\in [-\pi,\pi)^2$,
\begin{flalign}
   \left|   b_{\t}\big(c+cN_{\t},\tilde{c}\big) - (c,\tilde{c}) \right| & =\Big| b_\t(N_\t c,\tilde{c})  +|\t'|^2 |B|c\overline{\tilde{c}}\, \Big| \nonumber \\
   & \lesssim |\t| |c| |\tilde{c}| 
   .\label{24.10.26e10}
\end{flalign}
Then, \eqref{12.10.26e3} follows from \eqref{24.10.26e6}-\eqref{24.10.26e10} and \eqref{gatherkeystuff}.

It remains to notice that \eqref{9.10.26approx1}, \eqref{diagmm}, \eqref{12.10.26e3} and the definition of $\Phi$, prove Proposition \ref{prop.finalrelbd}.
\end{proof}

Let us finish this section with the main theorem which  collates  the previous results.
\begin{theorem}\label{thm.main1}
Let $u_{\ep,\t}$ solve \eqref{p0} and $v^{(2)}=c+z$ solve \eqref{prob.hom}. Then
\[
\| u_{\ep,\t}-v^{(2)} \|_{L^2(\square)} \lesssim \ep \| f\|_{L^2(\square)}.
\]
\end{theorem}
\begin{proof}
    This follows from the triangle inequality, \eqref{IH1}, \eqref{IH2} and \eqref{IH3}. 
\end{proof}

\section{Spectral analysis}\label{sec.spect}
\subsection{Comparison to a self-similar problem}
Consider the form
\begin{flalign}\label{Aform}
   A_{\xi}(c+z,\tilde{c}+\tilde{z})  :=  A^h \xi \cdot \xi  c \overline{\tilde{c}}    +  a^h|\xi_3|^2 (c+z,\tilde{c}+\tilde{z})_{L^2(B)}+   \big(\nabla'z,\nabla'\tilde{z}\big)_{L^2(B)}+
   (c+z,\tilde{c}+\tilde{z}),
\end{flalign}
on $V$,    and the unitary operator $\E_{\t'}:L^2(\square)\rightarrow L^2(\square)$ that is multiplication by the function $\exp(\i \t' \cdot y')$.
Clearly, one has
\begin{equation}\label{Ebdd}
\|   I - \E_{\t'} \|_{L^2(\square)\rightarrow L^2(\square)} \le |\t'|.  
\end{equation}

    Consider $v^{(3)} = c+z \in \CC\dot{+}H^1_0(B)$ the solution to 
\begin{equation}\label{prob.sshom}
A_{\t/\ep}(c+z,\tilde{c}+\tilde{z}) = (\E_{\t'}f,\tilde{c} +\tilde{z}), \quad \forall \tilde{c} \in \CC, \forall \tilde{z} \in H^1_0(B).   
\end{equation}
Clearly, one has \begin{equation}\label{ss.apr}
(1+|\xi|^2)|\tilde{c}|^2+ (1+|\xi_3|^2)\| \tilde{z}\|_{L^2(B)}^2 \lesssim A_\xi[\tilde{c} + \tilde{z}].
\end{equation}
\begin{lemma}\label{lem.ss}
Let $v^{(2)}$ solve \eqref{prob.hom} and $v^{(3)}$ solve \eqref{prob.sshom}. Then
\[
\| \E_{\t'}v^{(2)}-v^{(3)} \|_{L^2(\square)} \lesssim \ep \| f\|_{L^2(\square)}.
\]
\end{lemma}
\begin{proof}
Notice that $H^1_0(B)$ is an invariant subspace of $\E_{\t'}$, and we write  $v^{(2)}=c_2+ \E_{\t'}^*z_2$. Then, Proposition \ref{prop.homapriori}  implies
\begin{equation}\label{29.10.26e3}
\big(1 + \ep^{-2} |\t|^2 \big)|c_2| +  \big(1 + \ep^{-2} |\t_3|^2 \big)\|z_2\|_{L^2(B)} +\big(1 + \ep^{-1} |\t_3| \big)\|\nabla' z_2\|_{L^2(B)}  \lesssim \| f\|_{L^2(\square)}.
\end{equation}
Additionally, by the definition of the forms $A^h_{\ep,\t}$ and $A_\xi$, one has
\begin{multline*}
    A^h_{\ep,\t}(v^{(2)},\tilde{c}+\E_{\t'}^*\tilde{z})  \\ = A^h \frac{\t}{\ep} \cdot \frac{\t}{\ep} c_2 \overline{\tilde{c}} + a^h \left|\frac{\t_3}{\ep} \right|^2 (c_2 + \E_{\t'}^* z_2, \tilde{c} +\E_{\t'}^*\tilde{z})_{L^2(B)} + \big(\nabla' z_2, \nabla' \ \tilde{z})_{L^2(B)} +  (c_2 + \E_{\t'}^* z_2, \tilde{c} +\E_{\t'}^*\tilde{z}).
  \end{multline*}
Therefore, $c_2+z_2$ solves
\[
A_{\t/\ep}(c_2+z_2, \tilde{c} + \tilde{z})= (\E_{\t'}f,\tilde{c} +\tilde{z}) + R(\tilde{c}+\tilde{z}),
\]
for 
\begin{multline*}
     R(\tilde{c}+\tilde{z}) : = (f,\tilde{c} + \E_{\t'}^*\tilde{z}) - (\E_{\t'}f,\tilde{c} +\tilde{z}) +  (c_2+z_2,\tilde{c} + \tilde{z})  -(c_2+\E_{\t'}^*z_2,\tilde{c} + \E_{\t'}^*\tilde{z})  \\ + a^h\left|\tfrac{\t_3}{\ep}\right|^2 (c_2+z_2,\tilde{c}+\tilde{z})_{L^2(B)}  - a^h\left|\tfrac{\t_3}{\ep}\right|^2 (c_2+\E_{\t'}^*z_2,\tilde{c}+\E_{\t'}^* \tilde{z})_{L^2(B)}. 
\end{multline*}
Therefore, 
\begin{flalign}\label{29.10.26e1}
  A_{\t/\ep}[c_2+z_2-v^{(3)}] & =   R\big(c_2+z_2-v^{(3)}\big).
\end{flalign}
Now,
\begin{multline*}
 R(\tilde{c}+\tilde{z}) =\big((I-\E_{\t'})f,\tilde{c} \big)+ \big(c_2,(I-\E_{\t'}^*)\tilde{z} \big)+ \big((I-\E_{\t'}^*)z_2,\tilde{c} \big) \\+a^h\left|\tfrac{\t_3}{\ep}\right|^2 (c_2,(I-\E_{\t'}^*)\tilde{z})_{L^2(B)}+a^h\left|\tfrac{\t_3}{\ep}\right|^2 \big((I-\E_{\t'}^*)z_2,\tilde{c}\big)_{L^2(B)}    
\end{multline*}
so
\begin{flalign}
 &   \big|    R(\tilde{c}+\tilde{z})\big|\nonumber\\
        & \lesssim |\t'|  \| f\|_{L^2(\square)} |\tilde{c}| + |\t'|\big(  |c_2| \| \tilde{z}\|_{L^2(\square)} + \| z_2\|_{L^2(\square)} |\tilde{c}|  \big)  +\ep^{-2} |\t_3|^2|\t'| \big(  |c_2| \| \tilde{z}\|_{L^2(\square)} + \| z_2\|_{L^2(\square)} |\tilde{c}|  \big) \nonumber \\
        & \lesssim  \|f\|_{L^2(\square)} \big(|\t'| |\tilde{c}| + \ep \| \tilde{z}\|_{L^2(B)}  +|\t'| |\tilde{c}| + |\t_3|\| \tilde{z}\|_{L^2(B)} + |\t'| |\tilde{c}| \big) \nonumber \\
           & \lesssim \ep  \|f\|_{L^2(\square)}\sqrt{A_{\t/\ep}[\tilde{c}+\tilde{z}]},\label{29.10.26e2}
\end{flalign}
where the last two inequalities used \eqref{29.10.26e3} and \eqref{ss.apr} respectively. Therefore, \eqref{29.10.26e1} and \eqref{29.10.26e2} give
\[
A_{\t/ \ep}[c_2+z_2-v^{(3)}] \lesssim \ep^2 \| f\|_{L^2(\square)}^2.
\]
In particular,
\begin{equation}\label{SC1}
    \| c_2+z_2-v^{(3)}\|_{L^2(\square)} \lesssim \ep \| f\|_{L^2(\square)}.
\end{equation}
It remains to note, by \eqref{Ebdd} and \eqref{29.10.26e3}, that
\begin{equation}\label{SC2}
    \| \E_{\t'} c_2 - c_2\|_{L^2(\square)} \lesssim |\t'| | c_2| \lesssim \ep \| f\|_{L^2(\square)}.
\end{equation}
Then, the fact that $\E_{\t'} v^{(2)} = \E_{\t'} c_2 +z_2$, \eqref{SC1} and \eqref{SC2} establish the desired result.
\end{proof}
Combining Theorem \ref{thm.main1} and Lemma \ref{lem.ss} gives the following result:
\begin{theorem}\label{thm.main2}
Let $u_{\ep,\t}$ solve \eqref{p0} and $v^{(3)}$ solve \eqref{prob.sshom}. Then
\[
\| u_{\ep,\t}-\E_{\t'}^*v^{(3)} \|_{L^2(\square)} \lesssim \ep \| f\|_{L^2(\square)}.
\]
\end{theorem}
\subsection{Approximation of spectra}\label{sec.spectapp}
Here, we will interpret Theorem \ref{thm.main2} in operator-theoretic language. To this end, 
for each $\ep$ and $\t$, we introduce $L_{\ep,\t}$, the unbounded self-adjoint operator, in $L^2(\square)$, associated with the positive closed form $\ep^{-2} a_{\t} +b_\t$ with form domain $H^1_{per}(\square)$. Then, the solution $u_{\ep,\t}$ to \eqref{p0} is precisely $u_{\ep,\t} =  L_{\ep,\t}^{-1}f$. In a similar manner, we interpret $v^{(3)}$. We denote by $\H$, the closure of $V$ (cf. \eqref{V}) in $L^2(\square)$, that is
\[
\H = \{ h \in L^2(\square) \, | \, \text{$h = c + h_0$ where $c \in \CC$ and $h_0$ has support in $F$ and is independent of $y_3$} \}.
\]
Furthermore, we consider the orthogonal projector $\P:L^2(\square) \rightarrow \H$:
\[
\P u := c + h_0, \quad c= |\square\backslash F|^{-1} \int_{\square \backslash F}  u, \qquad  h_0(y)= \chi_0(y')  \Big(  \int_{-1/2}^{1/2} u(y',y_3) \, \mathrm{d}y_3 - c \Big), \ \text{for}\  y\in \square.
\]
Moreover, for each $\xi \in \RR^3$, let $\L_\xi$ be the unbounded self-adjoint operator, in $\H$,   generated by the positive form \eqref{Aform} with form domain $V$. Then, 
\[
v^{(3)}= \L_{\t/\ep}^{-1} \P \E_{\t'}f= \P^*\L_{\t/\ep}^{-1} \P \E_{\t'}f.
\]
Therefore, Theorem \ref{thm.main2} can be reformulated as follows:
\begin{theorem}\label{thm1.2}
	There exists a positive constant $C$ such that the inequality
	\[
	\| L_{\ep,\t}^{-1} - \E_{\t'}^*\P^*\L_{\t/\ep}^{-1} \P \E_{\t'} \|_{L^2(\square) \rightarrow L^2(\square)} \le C \ep
	\]
holds for all  $0<\ep<1$ and $\t\in \square^*$. 
\end{theorem}
Note that the operators $L_{\ep,\t}^{-1}$ and $\L_\xi^{-1}$ are compact, which follows from the compact embedding of each form domain into $L^2(\square)$ and $\H$ respectively. Consequently, Theorem \ref{thm1.2} readily implies uniform approximations for the eigenvalues  $\{\lambda_{\ep}^{(k)}(\t) \}_{k=1}^\infty$ and $\{\Lambda^{(k)}(\xi) \}_{k=1}^\infty$  of $L_{\ep,\t}$ and $\L_\xi$ respectively, labelled in increasing order and repeated according to multiplicity. Indeed, one has the following result: 
\begin{theorem}\label{thm1.3}
	There exists a constant $C>0$ such that 
\[
\left| \frac{1}{\lambda_{\ep}^{(k)}(\t)} - \frac{1}{\Lambda^{(k)}\left(\tfrac{\t}{\ep}\right)} \right| \le C \ep, \quad   \forall \t \in \square^*, \, \forall k \in \NN, \, 0<\ep<1.
\]
\end{theorem} 
\begin{proof}
By Theorem \ref{thm1.2}, the $k$-th eigenvalue of $L^{-1}_{\ep,\t}$ is $\ep$-close to the $k$-th eigenvalue of  the compact operator $ \E_{\t'}^*\P^*\L_{\xi}^{-1} \P \E_{\t'}$. It remains to notice that
    \begin{flalign*}
        \frac{\big( \E_{\t'}^*\P^*\L_{\xi}^{-1} \P \E_{\t'} f,f \big)}{\| f\|^2} =  \frac{( \L_{\xi}^{-1} \P \E_{\t'} f,\P \E_{\t'}f )}{\| \E_{\t'}f\|^2} =  \frac{( \L_{\xi}^{-1} \P \E_{\t'} f,\P \E_{\t'}f )}{\| \P \E_{\t'}f\|^2+ \| (I-\P) \E_{\t'} f\|_{L^2(\square)}^2},
    \end{flalign*}
and therefore, by the min-max principle and the fact that $\H$ is infinite dimensional, one can see that $k$-th eigenvalue of $\E_{\t'}^*\P^*\L_{\xi}^{-1} \P \E_{\t'} $ coincides with the $k$-th eigenvalue of $\L_{\xi}^{-1}$ for all $k\in \NN$.
\end{proof}
We conclude by noting that, following the standard procedure (see, e.g., \cite{Zhi2005,CoKaSm}), one can show that, if the eigenvalues, $\mu_n$, of the Dirichlet Laplacian on $B$ are simple and the corresponding normalised eigenfunctions, $\phi_n$, have non-zero mean,  then the eigenvalues $\Lambda^{(k)}(\xi)$ are solutions $\lambda$ of the dispersion relation
\[
A^h\xi\cdot\xi-(1-|B|)\,a^h|\xi_3|^2=\beta\big(\lambda-1-a^h|\xi_3|^2\big),
\]
where
\[
\beta(\mu):=\mu+\mu^2\sum_{n=1}^\infty\frac{|(\phi_n,1)_{L^2(B)}|^2}{\mu_n-\mu}.
\]

\section{Comparison to the two-scale operator}\label{sec.twoscale}
In this section, we show that $\L_\xi$ are the fibres of the two-scale operator, associated with problem \eqref{p4} after performing Fourier transform in the macroscopic variable $x$. We will then use this observation, in conjunction with Theorem \ref{thm1.2},  to approximate the solution of the original problem \eqref{p.strong} in terms of the two-scale homogenized limit problem \eqref{p4}. 

To this end, let us introduce the Hilbert space 
\[
L^2(\RR^3;\mathcal{H}),
\]
$P : L^2(\RR^3\times \square) \rightarrow L^2(\RR^3;\mathcal{H})$ the orthogonal projection, defined for example by the continuous extension of
\[
P(f(x)g(y)) = f(x)\P g(y),
\]
and define the dense subspace
\begin{multline*}
\mathbb{V} : = \{ u+v \, : \text{for  $ u \in H^1(\RR^3)$ and   $v \in L^2(\RR^3; H^1_{per}(\square))$ with $\partial_{x_3} v \in L^2(\RR^3 \times \square)$ and}\\ \text{ with  $v$  independent of $y_3$ and supported on $\RR^3 \times F$.} \}    
\end{multline*}
Finally, let $\mathbb{L}$ be the unbounded self-adjoint operator defined by the closed quadratic form
\begin{equation*}
	 	\begin{aligned}
& \mathbb{A}[u+v] = \int_{\RR^3} \Big(	 A^{\mathrm h} \nabla_x u \cdot \overline{\nabla_x u}  + \  a^{\mathrm h}  \|\partial_{x_3}( u+  v)\|^2_{L^2(B)}   +  \| \nabla_{y'} v\|_{L^2(B)}^2 + \| u + v\|_{L^2(\square)}^2 \Big)\, \mathrm{d}x,
	\end{aligned} 
\end{equation*}
with form domain $\mathbb{V}$.

As $\mathbb{L}$ is a constant coefficient differential operator in $x$, it admits a fibre decomposition with respect to the Fourier transform 
in the first variable. Namely, introducing the Fourier transform $\mathcal{F}: L^2(\RR^3 \times \square) \rightarrow L^2(\RR^3 \times \square)$
\[
\mathcal{F} f(\xi,y) = (2\pi)^{-3/2} \int_{\RR^3} f(x,y) e^{- \i x \cdot \xi} \, \mathrm{d} x;
\]
one has
\[
\mathcal{F}\,\mathbb{L}\,\mathcal{F}^{-1} = \int^{\oplus}_{\RR^3} \mathbb{L}_\xi \, \mathrm{d} \xi,
\]
where $\mathbb{L}_\xi$ is the unbounded linear operator generated by the form 
\[
 A^{\mathrm h} \xi \cdot \xi | c|^2   + \int_B   a^{\mathrm h} |\xi_3|^2 | c+  z(y')|^2 \, \mathrm{d}y'   +  \int_B  | \nabla_{y'} z(y')|^2 \, \mathrm{d}y' +\int_\square | c+  z(y')|^2 \, \mathrm{d}y,
\]
with form domain 
\[
\{ c+z \, : \text{for  $ c \in \CC$ and   $z \in H^1_{per}(\square)$ with $z$ independent of $y_3$ and supported on $F$.} \}  
\]
Therefore, we observe that $\mathbb{L}_\xi$ is exactly the operator $\L_\xi$ introduced in Section \ref{sec.spectapp}; consequently,
\[
\mathcal{F}\,\mathbb{L}^{-1}\,\mathcal{F}^{-1} = \int^{\oplus}_{\RR^3} \L_\xi^{-1} \, \mathrm{d} \xi.
\]
Moreover, one has
\[
P \mathcal{F}^{-1} = \mathcal{F}^{-1}\int_{\RR^3}^{\oplus} \P \, \mathrm{d}\xi, 
\]
and so
\[
\mathcal{F}P^* \mathbb{L}^{-1} P\,\mathcal{F}^{-1} = \int^{\oplus}_{\RR^3} \P^* \L_{\xi}^{-1} \P \, \mathrm{d} \xi.
\]
Now, for the unitary rescaling in the first variable, $\Gamma_\ep : L^2(\RR^3\times \square) \rightarrow L^2(\RR^3\times \square)$, $\Gamma_\ep f(\xi,y) = \ep^{3/2} f(\ep \xi,y)$, notice that
\[
\mathcal{F} \Gamma_\ep = \Gamma_{\ep}^{-1} \mathcal{F};
\]
and, therefore,
\[
\mathcal{F} \Gamma_\ep P^* \mathbb{L}^{-1} P\Gamma_\ep^{-1}\mathcal{F}^{-1} = \int^{\oplus}_{\RR^3} \P^* \L_{\xi/\ep}^{-1} \P \, \mathrm{d} \xi.
\]
Introducing the bounded linear operator $\chi : L^2(\square^* \times \square) \rightarrow L^2(\RR^3 \times \square)$ that extends functions by zero outside of $\square^*$; one has
 \[
\chi^* \mathcal{F} \Gamma_\ep P^* \mathbb{L}^{-1} P\Gamma_\ep^{-1}\mathcal{F}^{-1} \chi= \int^{\oplus}_{\square^*} \P^* \L_{\theta/\ep}^{-1} \P  \, \mathrm{d} \theta.
\]
Defining 
\[
\E =\int_{\square^*}^{\oplus} \E_{\t'} \, \mathrm{d} \t,
\]
gives
 \[
\E^*\chi^* \mathcal{F} \Gamma_\ep P^* \mathbb{L}^{-1} P\Gamma_\ep^{-1}\mathcal{F}^{-1} \chi \E= \int^{\oplus}_{\square^*} \E_{\t'}^* \P^* \L_{\theta/\ep}^{-1} \P \E_{\t'}  \, \mathrm{d} \theta.
\]
Then, Theorem \ref{thm1.2} implies
\[
	\|\int_{\square^*}^{\oplus} L_{\ep,\t}^{-1} \, \mathrm{d}\theta-  \E^*\chi^* \mathcal{F} \Gamma_\ep P^* \mathbb{L}^{-1} P\Gamma_\ep^{-1}\mathcal{F}^{-1} \chi \E\|_{L^2(\square^*\times \square) \rightarrow L^2(\square^*\times \square)} \le C \ep.
	\]
    Recalling that, for the unitary rescaling $\Gamma_\ep$ and unitary Gelfand transform $U$, one has
    \[
U  \Gamma_\ep  L_\ep^{-1} \Gamma_{\ep}^{-1} U^{-1}  =   \int_{\square^*}^{\oplus} L_{\ep,\t}^{-1} \, \mathrm{d}\theta;
    \]
then upon introducing
    \[
J_\ep := \Gamma_\ep^{-1}\mathcal{F}^{-1}\chi \E U \Gamma_\ep
    \]    
    gives the following result: 
      \begin{theorem}\label{thm2}
      For the constant $C>0$ in Theorem \ref{thm1.2},  the inequality
 	\[
 	\| L_\ep^{-1} - J_\ep^* P^* \mathbb{L}^{-1} P  J_\ep \|_{L^2(\RR^3) \rightarrow L^2(\RR^3)} \le C \ep
 	\]
 	holds for all $0<\ep<1$. 
 \end{theorem}
 \begin{remark}
     Since $\mathbb{L}^{-1}$ acts from $L^2(\RR^3;\H)$ to $L^2(\RR^3;\H)$, $P^*$ in the above assertions can be removed. 
 \end{remark}
 \begin{remark}
     One can readily show that
     \[
J_\ep^* J_\ep = I 
     \]
     and
     \[
J_\ep J_\ep^* = \Gamma_\ep^{-1} \mathcal{F}^{-1} \chi \chi^* \mathcal{F} \Gamma_\ep = \mathcal{F}^{-1} \Gamma_\ep \chi \chi^*\Gamma_\ep^{-1} \mathcal{F} = \mathcal{F}^{-1} \mathds{1}_{\xi \in \ep^{-1} \square^*} \mathcal{F}.
     \]
     Clearly $J_\ep J_\ep^*$ converges to $I$, as $\ep \rightarrow 0$, in the strong operator topology.
 \end{remark}

 \end{document}